\documentclass[a4paper, 11pt]{article}

\usepackage{graphicx} 
\usepackage[all]{xy} 
\usepackage{amsfonts, amssymb, amsmath, url} 
\usepackage{amscd}
\usepackage{stmaryrd}
\usepackage{amsthm}
\usepackage{ifthen}
\usepackage{epsfig}

\newcommand{\comment}[1]{}

\comment{
\usepackage{graphicx} 
\usepackage[all]{xy} 
\usepackage{amsfonts, amssymb, amsmath} 
\usepackage{amscd}
\usepackage{amsthm}
\usepackage{epsfig}
\usepackage{ifthen}
}

\newtheorem{Thm}{Theorem}[section]
\newtheorem{Lem}{Lemma}[section]
\newtheorem{Prop}{Proposition}[section]

\newtheorem{Question}{Question}[section]
\newtheorem{Coro}{Corollary}[section]

\theoremstyle{definition}
\newtheorem{Rem}{Remark}[section]
\newtheorem{Def}{Definition}[section]

\newtheorem{Example}{Example}[section]

\numberwithin{equation}{section}

\newcommand{\ind}{{\bf 1}}
\newcommand{\proba}{\,\mathbb P}
\newcommand{\esp}{\,{\mathbb E}}
\newcommand{\supp}{{\rm{supp}}}

\newcommand{\defe}{\mathrel{\mathop:}=}

\newcommand{\defd}{\eqd}

\newcommand{\inv}{^{-1}}
\newcommand{\var}{{\rm{Var}}}

\newcommand{\calA}{{\cal A}}
\newcommand{\calB}{{\cal B}}

\newcommand{\calE}{{\cal E}}
\newcommand{\filF}{{\cal F}}
\newcommand{\calG}{{\cal G}}

\newcommand{\calP}{{\cal P}}

\newcommand{\calR}{{\cal R}}
\newcommand{\calS}{{\cal S}}
\newcommand{\calT}{{\cal T}}

\def\indt#1{\{#1_t\}_{t\in T}}
\def\indtr#1{\{#1_t\}_{t\in\mathbb R}}

\def\indtauto#1{\left\{#1_t\right\}_{t\in T}}

\def\indn#1{\{#1_n\}_{n\in \mathbb N}}

\newcommand{\msspan}{\vee\mbox{-}{\rm{span}}}
\newcommand{\cmsspan}{\overline{\vee\mbox{-}{\rm{span}}}}

\def\laps{{L^\alpha_+(S,\mu)}}

\def\lapdo{{L^\alpha_+(S_1,\mu_1)}}
\def\lapim{{L^\alpha_+(S_2,\mu_2)}}
\def\lapdoSC{L^\alpha_+(S_1,\sigma(C),\mu_1)}
\def\lapimSC{L^\alpha_+(S_2,\sigma(U(C)),\mu_2)}
\def\lapdoS{L^\alpha_+(S_1,\sigma(\filF),\mu_1)}
\def\lapimS{L^\alpha_+(S_2,\sigma(U(\filF)),\mu_2)}
\def\lapdoN{{L^\alpha_+(\mu_1)}}
\def\lapimN{{L^\alpha_+(\mu_2)}}
\def\inddo{\ind_{S_1}}
\def\indim{\ind_{S_2}}

\def\lap{L^\alpha_+}
\def\la{L^\alpha}

\newcommand{\eqnh}{\begin{eqnarray*}}
\newcommand{\eqne}{\end{eqnarray*}}
\newcommand{\eqnhn}{\begin{eqnarray}}
\newcommand{\eqnen}{\end{eqnarray}}
\newcommand{\equh}{\begin{equation}}
\newcommand{\eque}{\end{equation}}

\newcommand{\sumin}{\sum_{i=1}^n}

\newcommand{\eqd}{\stackrel{\rm d}{=}}

\newcommand{\widebar}{\overline}
\newcommand{\eintt}{\ \int^{\!\!\!\!\!\!\!e}}
\newcommand{\Eintt}{\int^{\!\!\!\!\!\!\!e}}
\def\eint#1{\, \int^{\!\!\!\!\!\!\!e}_{#1}}
\def\Eint#1{\, \int^{\!\!\!\!\!\!\!\!{e}}_{#1}}

\def\topp#1{^{(#1)}}
\def\dfrac#1#2#3{\frac{d(#1\circ #2)}{d#3}}

\def\lap{L^\alpha_+}

\def\lopT{L^0_+(T,\calT,\lambda)}

\def\tia{^{1/\alpha}}
\def\ae{\mbox{-a.e.}}

\def\malmosts{\mbox{ almost surely}}
\def\inftydots#1{#1_1,#1_2,\dots}

\def\mrn{\mathbb R^{\mathbb N}}
\def\mrho{\rho_{\mu,\alpha}}

\def\ratios#1{\calR_+({#1})}
\def\ratiosf{\calR_+(\filF)}
\def\ratiosuf{\calR_+(U(\filF))}
\def\eratios#1{\calR_{e,+}({#1})}
\def\eratiosf{\calR_{e,+}(\filF)}
\def\eratiosuf{\calR_{e,+}(U(\filF))}
\def\ssS{{S_{I,N}}}
\def\ssL{{\lambda_{I,N}}}

\def\ssB{{\calB_\ssS}}
\def\ssF{{f_{I,N}}}
\def\xI{X^{I}}
\def\xN{X^{N}}

\def\bccbb#1{\Big\{#1\Big\}}

\def\tpd#1{^{\rm #1}}

\def\rnp{{\mathbb R_+^{\mathbb N}}}
\def\d{{\rm d}}
\def\e{{\rm e}}

\def\eqnhspace{& & \ \ \ \ }

\def\noi{\noindent}
\def\calF{{\mathcal F}}
\def\bbR{{\mathbb R}}
\def\bbQ{{\mathbb Q}}
\def\bbN{{\mathbb N}}
\def\P{{\mathbb P}}

\def\mand{\mbox{ and }}
\def\qmand{\quad\mbox{ and }\quad}

\def\smu{{(S,\mu)}}

\def\adaptF#1{\{#1_t,\filF_t:0\leq t<\infty\}}

\newcommand{\pxt}{\{X_t\}_{t\in T}}
\def\itemnumber#1{\noindent\parbox{0.27in}{\it (#1)}}
\newcommand{\fft}{\{f_t\}_{t\in T}}

\def\wtilde{\widetilde}
\def\E{{\mathbb E}}

\def\noi{\noindent}
\def\calF{{\mathcal F}}
\def\bbR{{\mathbb R}}
\def\bbQ{{\mathbb Q}}
\def\bbN{{\mathbb N}}
\def\bbZ{{\mathbb Z}}
\def\P{{\mathbb P}}

\newboolean{qedTrue}
\setboolean{qedTrue}{false}

\begin{document}
\title{{On the Structure and Representations of Max-Stable Processes}\thanks{The authors were partially supported by NSF grant DMS--0806094 at the University of Michigan.}}
\author{Yizao Wang and Stilian A. Stoev
\thanks{{\it Address:}
 Department of Statistics, The University of Michigan, 439 W.\ Hall, 1085 S.\ University, Ann Arbor,
  MI 48109--1107; 
  {E--mails:} 
  {\texttt{\{yizwang, sstoev\}}}@{\texttt umich.edu}.
  } 
}
\maketitle

\comment{
\title{On the Structure and Representations of Max--Stable Processes\protect\thanksref{T1}}
\runtitle{Structure of Max--Stable Processes}


\author{\fnms{Yizao} \snm{Wang}\ead[label=e1]{yizwang@umich.edu}},
\and
\author{\fnms{Stilian A.} \snm{Stoev}\ead[label=e2]{sstoev@umich.edu}}

\thankstext{T1}{The authors were partially supported by NSF grant DMS--0806094 at the University of Michigan.}
\runauthor{Y. Wang and S. Stoev}

\affiliation{Department of Statistics, University of Michigan}

\address{ Department of Statistics,\\
The University of Michigan, \\
439 W.\ Hall, 1085 S.\ University, \\
Ann Arbor, MI 48109--1107\\
\printead{e1}\\
\phantom{E-mail:\ }\printead*{e2}
}
}
\begin{abstract}
We develop classification results for max--stable processes, based on their spectral representations.  The structure of max--linear isometries and minimal spectral 
representations play important roles. 
We propose a general classification strategy for measurable max--stable processes based on the notion of co--spectral functions.  
In particular, we discuss the spectrally continuous--discrete, the conservative--dissipative, and positive--null decompositions.
For stationary max--stable processes, the latter two decompositions arise from connections to non--singular flows and are closely related 
to the classification of stationary sum--stable processes. The interplay between the introduced decompositions of max--stable processes is further explored.
As an example, the Brown--Resnick stationary processes, driven by fractional Brownian motions, are shown to be dissipative.  A result on general Gaussian processes with 
stationary increments and continuous paths is obtained.
\end{abstract}

\comment{
\begin{keyword}[class=AMS]
\kwd[Primary ]{60G52}
\kwd{60G70}
\kwd[; secondary ]{37A50}
\end{keyword}

\begin{keyword}
\kwd{max--stable}
\kwd{classification}
\kwd{max--linear isometry}
\kwd{spectral representation}
\kwd{non--singular flow}
\kwd{co--spectral function}
\end{keyword}
}


\section{Introduction}

Max--stable processes have been studied extensively in the past 30
years.  The works of Balkema and Resnick~\cite{balkema77max}, de
Haan~\cite{dehaan78characterization,dehaan84spectral}, de Haan and
Pickands~\cite{dehaan86stationary}, Gin\'e {\it et al.}~\cite{gine90max} and
Resnick and Roy~\cite{resnick91random}, among many others have lead to
a wealth of knowledge on max--stable processes.
The seminal works of de Haan \cite{dehaan84spectral} and de Haan and
Pickands \cite{dehaan86stationary} laid the foundations of the
spectral representations of max--stable processes and established
important structural results for stationary max--stable processes.
Since then, however, while many authors focused on various important
aspects of max--stable processes, the general theory of their
representation and structural properties had not been thoroughly
explored.  At the same time, the structure and the classification of
sum--stable processes has been vigorously studied.  Rosi\'nski \cite{rosinski95structure}, building
on the seminal works of Hardin \cite{hardin81isometries,hardin82spectral} about minimal
representations, developed the important connection between stationary
sum--stable processes and flows.  This lead to a number of
important contributions on the structure of sum--stable processes (see, e.g.\
\cite{rosinski96classes,rosinski00decomposition,pipiras02structure,pipiras04stable,samorodnitsky05null}).
There are relatively few results of this nature about the structure of max--stable processes, with 
the notable exceptions of de Haan and Pickands \cite{dehaan86stationary}, Davis and Resnick \cite{davis93prediction}
and the very recent works of Kabluchko {\it et al.}~\cite{kabluchko08stationary} and Kabluchko~\cite{kabluchko08spectral}.

Our goal here is to develop representation and classification theory
for max--stable processes, similar to the available one for sum--stable processes.
We are motivated by the strong similarities between the spectral representations of 
sum-- and max--stable processes.  This procedure however, is non--trivial.  
The notion of {\it minimal extremal integral} representation plays a key role
as does the {\it minimal integral representation} for $\alpha$--stable processes (see
Hardin~\cite{hardin82spectral} and Rosi\'nski \cite{rosinski95structure,rosinski06minimal}). 
Before one can fruitfully handle the {\it minimal extremal integral} representations,
it turns out that one should first thoroughly investigate the
structure of max--linear isometries, also known as the {\it pistons} of
de Haan and Pickands~\cite{dehaan86stationary}.  We refine and extend their work in Section~\ref{sec:maxLinear}.  In Section \ref{sec:minimal}, we develop the theory of minimal
representations for max--stable processes.  Our approach is motivated by the works of Hardin \cite{hardin82spectral} and
Rosi\'nski \cite{rosinski95structure} in the sum--stable context. 

In Section \ref{sec:classification}, we establish general classification results for max--stable processes
by using the developed theory of {\it minimal spectral representations}. In Section~\ref{sec:continuousDiscreteDecomposition}, we first show that essentially any
max--stable process can be represented uniquely as the maximum of two independent components, characterized as
{\it spectrally continuous} and {\it spectrally discrete}, respectively. The spectrally discrete part gives rise 
to the notion of {\it discrete principal components}, which may be of independent interest in modeling of max--stable
processes and fields.  

In Section~\ref{sec:cospectral}, we introduce the notion of {\it co--spectral functions}, for the large class of measurable max--stable processes $X = \indt X$. There $T$ is a separable metric space equipped with the Borel--$\sigma$--algebra and a $\sigma$--finite measure. The co--spectral functions of such processes are invariant to the choice of the spectral representations, up to a multiplicative factor. This allows us to develop a general strategy for the classification of measurable $\alpha$--Fr\'echet processes, based on positive cones of co--spectral functions. As particular examples, we obtain the {\it conservative--dissipative} and
{\it positive--null} decompositions, which correspond to certain choices of cones for the 
co--spectral functions.

Section \ref{sec:stationary} is devoted to the classification of stationary max--stable processes.
As in the sum--stable case, the minimal representations allow us to associate a measurable 
non--singular flow to every measurable stationary max--stable process.  This correspondence enables one to
apply existing ergodic theory results about the flow to characterize the max--stable process.
The conservative--dissipative and positive--null decompositions introduced in Examples 
\ref{sec:cons-diss} and \ref{sec:pos-null} are in fact motivated by the corresponding
decompositions of the underlying flow. These two results are in close correspondence with the
classifications of Rosi\'nski \cite{rosinski95structure} and Samorodnitsky
\cite{samorodnitsky05null} for sum--stable processes.  As in Rosi\'nski \cite{rosinski95structure},
we obtain that the class of stationary max--stable processes generated by dissipative flows is precisely
the class of mixed moving maxima.  

In Section~\ref{sec:BRp}, we apply the results in Section~\ref{sec:stationary} to Brown--Resnick processes. We give simple necessary and sufficient conditions 
for a generalized Brown--Resnick stationary process to be a mixed moving maxima.  This extends and complements 
the recent results of Kabluchko {\it et al.}~\cite{kabluchko08stationary}.  In fact, as a by--product, 
by combining our results and those in \cite{kabluchko08stationary}, we obtain an interesting
fact about general zero--mean Gaussian processes $W=\{W_t\}_{t\in \bbR}$ with stationary increments and continuous paths.  Namely, 
for such processes, we have that, with probability one,
$$
\lim_{|t|\to\infty} {\Big(} W_t -{\rm Var}(W_t)/2 {\Big)}  = -\infty \ \mbox{ {\it implies} }\ 
 \int_{\bbR} \exp\{ W_t - {\rm Var}(W_t)/2\} dt < \infty.
$$
In particular, we show that if $\indtr W$ is a fractional Brownian motion, then the generated Brown--Resnick process is a mixed moving maxima.
We conclude Section~\ref{sec:BRp} with some open questions.  Some proofs and auxiliary results are given in
the Appendix.

Part of our results in Sections~\ref{sec:classification} and~\ref{sec:stationary} are modifications and extensions of results of de Haan and Pickands~\cite{dehaan86stationary}. The main difference is that we provide a complete treatment of the measurability issue, when the processes are continuously indexed. 
Before we proceed with the more technical preliminaries, we are obliged to mention the recent work of
Kabluchko \cite{kabluchko08spectral}.  In this exciting contribution, the author establishes some very 
similar classification results by using an {\it association device} between max-- and sum--stable processes.
This association allows one to transfer existing classifications of sum--stable processes to the max--stable domain.
It also clarifies the connection between these two classes of processes.  Our results were obtained 
independently and by using rather different technical tools.  The combination of the two approaches provides a more clear 
picture on the structure of max-- and sum--stable processes as well as their interplay.

\section{Preliminaries}\label{sec:prelim}

The importance of max--stable processes stems from the fact that they
arise in the limit of the component--wise maxima of independent
processes.  It is well known that the univariate marginals of a 
max--stable process are necessarily extreme value distributions, i.e.\ 
up to rescaling and shift they are either Fr\'echet, Gumbel or negative Fr\'echet.
The dependence structure of the max--stable processes, however, can be quite intricate and it does
not hinge on the extreme value type of the marginal distributions (see e.g.\ Proposition 5.11 in 
Resnick \cite{resnick87extreme}).  Therefore, for convenience and without loss of generality we will
focus here on max--stable process with Fr\'echet marginal
distributions.  Recall that a positive random variable $Z\ge 0$ has
$\alpha$--Fr\'echet distribution, $\alpha>0$, if
\[
\proba(Z\leq x) = \exp\{-\sigma^\alpha x^{-\alpha}\}\,,x\in(0,\infty)\,.
\]
Here $\left\| Z \right\|_\alpha \defe \sigma>0$ stands
for the \textit{scale coefficient} of $Z$.
It turns out that a stochastic process $\indt X$ with
$\alpha$--Fr\'echet marginals is max--stable {\it if and only if} all
positive {\it max--linear combinations}:
\begin{equation}\label{e:max-lin}
  \max_{1\leq j\leq n}a_jX_{t_j} \equiv \bigvee_{1\leq j\leq n}a_jX_{t_j}\ \ \ \forall a_j> 0,\ t_j\in T,\ 1\le j\le n,
\end{equation}
are $\alpha$--Fr\'echet random variables (see de
Haan~\cite{dehaan78characterization} and e.g.\ \cite{stoev06extremal}).  This feature resembles the definition of
Gaussian or, more generally, symmetric $\alpha$--stable (sum--stable)
processes, where all finite--dimensional linear combinations are
univariate Gaussian or symmetric $\alpha$--stable, respectively (see e.g.\ \cite{samorodnitsky94stable}).  We
shall therefore refer to the max--stable processes with $\alpha$--Fr\'echet marginals as to {\it $\alpha$--Fr\'echet processes}.

The seminal work of de Haan \cite{dehaan84spectral} provides convenient {\it spectral representations}
for stochastically continuous $\alpha$--Fr\'echet processes in terms of functionals of 
Poisson point processes on $(0,1)\times (0,\infty)$.  Here, we adopt the slightly more general, but
essentially equivalent, approach of representing max--stable processes through extremal integrals 
with respect to a random sup--measures (see Stoev and Taqqu \cite{stoev06extremal}). We do so in order to emphasize the analogies with the well--developed
theory of sum--stable processes (see e.g.\ Samorodnitsky and Taqqu \cite{samorodnitsky94stable}). 

\begin{Def} \label{def:M_alpha} Consider a measure space $(S,{\cal S},\mu)$ and suppose $\alpha>0$.  A stochastic process $\{M_\alpha(A)\}_{A\in {\cal S}}$,
indexed by the measurable sets $A\in {\cal S}$ is said to be an \textit{$\alpha$--Fr\'echet random sup--measure} with \textit{control measure} $\mu$, if the following 
conditions hold:\\
\itemnumber  i the $M_\alpha(A_i)$'s are independent for disjoint $A_i\in {\cal S},\ 1\le i\le n.$\\
\itemnumber {ii} $M_\alpha(A)$ is $\alpha$--Fr\'echet with scale coefficient $\|M_\alpha(A)\|_\alpha= \mu(A)^{1/\alpha}$.\\
\itemnumber{iii} for all disjoint $A_i$'s, $i\in \bbN$, we have
 $M_\alpha(\cup_{i\in \bbN} A_i) = \bigvee_{i\in \bbN} M_\alpha(A_i),$ almost surely.
\end{Def}

\noindent
Now, given an $\alpha$--Fr\'echet random sup--measure $M_\alpha$ as above, one can define the {\it extremal integral}
of a non--negative simple function  $f(u):= \sum_{i=1}^n a_i 1_{A_i}(u) \ge 0,\ A_i\in {\cal S}$:
$$
 \Eint{S} fdM_\alpha \equiv \Eint{S} f(u) M_\alpha(du) := \bigvee_{1\le i\le n} a_i M_\alpha(A_i).
$$ 
The resulting extremal integral is an $\alpha$--Fr\'echet random variable with scale coefficient $(\int_E f^\alpha
d\mu)^{1/\alpha}$. The definition of $\eint{S}fdM_\alpha$ can, by continuity in probability, be naturally 
extended to integrands $f$ in the space
$$
 L_+^\alpha(S,\mu):={\Big\{}f:S\to \bbR_+\, :\, \mbox{ $f$ measurable with }\int_S f^\alpha d\mu <\infty{\Big\}}.
$$
It turns out that the random variables $\xi_j:=\eint{S}f_j dM_\alpha,\ 1\le j\le n$ are independent if and only if
the $f_j$'s have pairwise disjoint supports (mod $\mu$). Furthermore, the extremal integral is {\it  max--linear}:
$$
 \Eint{S} (a f\vee b g) d M_\alpha = a \Eint{S}  f dM_\alpha \vee  b  \Eint{S} g d M_\alpha,
$$
for all $a,b>0$ and $f,g \in L_+^\alpha(S,\mu).$ For more details, see Stoev and Taqqu \cite{stoev06extremal}.

Now, for any collection of deterministic functions $\indt f\subset\lap(S,\mu)$, one can construct the stochastic process: 
\equh\label{eq:xtsft}
X_t = \eintt_Sf_t(u)M_\alpha(du)\,,\forall t\in T\,.
\eque
In view of the max--linearity of the extremal integrals and \eqref{e:max-lin}, the resulting process $X=\{X_t\}_{t\in T}$ is $\alpha$--Fr\'echet.
Furthermore, for any $n\in\mathbb N,\ x_i>0,\ t_i\in T,\ 1\le i \le n$:
\equh\label{eq:xt1a1}
\proba\{ X_{t_1}\leq x_1,\dots,X_{t_n}\leq x_n \} = \exp{\Big\{}-\int_S\Big(\vee_{1\leq i\leq n} x_i\inv f_{t_i}(u)\Big)^\alpha \mu(du) {\Big\}}.
\eque
This shows that the deterministic functions $\indt f$ characterize completely the finite--dimensional distributions of the process $\indt X$.  In general, if
\begin{equation}\label{rep:extremalRep}
 \{X_t\}_{t\in T} \eqd {\Big\{}\Eint{S} f_t dM_\alpha {\Big\}}_{t\in T},
\end{equation}
for some $\indt f \subset L_+^\alpha(S,\mu)$, we shall say that the process $X=\{X_t\}_{t\in T}$ has the {\it extremal integral} or {\it spectral representation}
$\{f_t\}_{t\in T}$ over the space $L_+^\alpha(S,\mu)$.  The $f_t$'s in \eqref{rep:extremalRep} are also referred to as \textit{spectral functions} of $X$.

{\it Our goal in this paper is to characterize $\alpha$--Fr\'echet processes in terms of their spectral representations.}
Many $\alpha$--Fr\'echet processes of practical interest have tractable spectral representations. As shown in the proposition below, an $\alpha$--Fr\'echet process $X$ has the representation  \eqref{rep:extremalRep}, where $(S,\mu)$ is a {\it standard Lebesgue space} (see Appendix A in~\cite{pipiras04stable}), if and only if, $X$ satisfies {\it Condition S}.

\begin{Def}\label{d:Cond-S} An $\alpha$--Fr\'echet process $X = \{X_t\}_{t\in T}$ is said to satisfy {\it Condition S} if there exists a
countable subset $T_0\subseteq T$ such that for every $t\in T$, we have that 
$X_{t_n}\stackrel{P}\to X_t$ for some $\indn t\subset T_0$.
\end{Def}

\begin{Prop}\label{prop:conditionS} 
An $\alpha$--Fr\'echet process $X=\{X_t\}_{t\in T}$ has the extremal integral representation \eqref{rep:extremalRep}, with any (some) standard Lebesgue space $(S,\mu)$ and an $\alpha$--Fr\'echet random sup--measure on $S$ with control measure $\mu$,
if (only if) it satisfies Condition S.
\end{Prop}

\noindent The result above follows from Proposition 3.2 in \cite{stoev06extremal}, since the standard 
Lebesgue space $(S,\mu)$ may be chosen to be $[0,1]$, equipped with the Lebesgue
measure.
\begin{Rem}
As shown in Kabluchko \cite{kabluchko08spectral} (Theorem 1), every max--stable process can have a 
spectral representation over a sufficiently rich abstract measure space.
\end{Rem}

In the sequel, we focus only on the rich class of $\alpha$--Fr\'echet processes that satisfy 
Condition S.  This includes, for example, all measurable max--stable processes 
$X=\{X_t\}_{t\in T}$, indexed by a separable metric space $T$ 
(see Proposition \ref{p:measurability} below).

The fact that $(S,\mu)$ is a standard Lebesgue space implies that 
the space of integrands $L_+^\alpha(S,\mu)$ is a complete and {\it separable} metric space 
with respect to the metric:
\equh\label{eq:metric} 
 \mrho(f,g) = \int_S|f^\alpha-g^\alpha|d\mu\,.  
\eque
This metric is natural to use when handling extremal integrals, since as $n\to\infty$,
\begin{equation}\label{e:mrho}
\eint{S} f_n dM_\alpha \stackrel{P} {\longrightarrow} \xi\,,\ \ 
\mbox{{\it if and only if, }}\ \  \mrho(f_n,f) = \int_S |f_n^\alpha - f^\alpha| d\mu \to 0,\,
\end{equation}
where  $\xi = \Eint{S} f d M_\alpha$ (see e.g.\ \cite{stoev06extremal} and also 
Davis and Resnick \cite{davis93prediction}). 
In the sequel, we equip the space $L_+^\alpha(S,\mu)$ with the metric $\mrho$ and 
often write $\|f\|_{L_+^\alpha(S,\mu)}^\alpha$ for $\int_S f^\alpha d\mu$.

\section{Max--Linear Isometries}\label{sec:maxLinear}
The max--linear (sub)spaces of functions in $L_+^\alpha(S,\mu)$ play a key role in the
representation and characterization of max--stable processes.  We say that $\filF$ is a 
\textit{max--linear sub--space} of $\laps$ if the following conditions hold:\\
\itemnumber {i} $af\vee bg\in\filF$, for all $a,b>0,f,g\in\filF$.

\itemnumber {ii} $\filF\subset \laps$ is closed in the metric $\mrho$. \\
\noindent
In particular, we will frequently encounter the max--linear space $\filF \defe \cmsspan(f_t,t\in T)$, which is generated by 
the max--linear combinations $\vee_{1\le i\le n}a_i f_{t_i}$, $t_i\in T,\ a_i>0$, of the spectral 
functions in \eqref{rep:extremalRep}.  In view of \eqref{e:mrho}, the set of 
extremal integrals $\{\eint{S} f dM_\alpha,\ f\in {\cal F}\}$ is the smallest set that is closed 
with respect to convergence in probability and contains all max--linear combinations
$\vee_{1\le i\le n} a_i X_{t_i}$.  For more details, see \cite{stoev06extremal}.

An $\alpha$--Fr\'echet process $X =\{X_t\}_{t\in T}$ as in \eqref{eq:xtsft} has many equivalent 
spectral representations.  They are all related, however, through {\it max--linear isometries} (see 
e.g.\ \eqref{eq:Ucanonical} below):

\begin{Def} Let $\alpha>0$.  The map $U:\lap(S_1,\mu_1)\to \lap(S_2,\mu_2)$, is said to be 
a max--linear isometry, if:\\
\itemnumber i  $U(a_1f_1\vee a_2\,f_1) = a_1(Uf_1)\vee a_2(Uf_2), \mu_2\ae$, for all $f_1,f_2\in \lap(S_1,\mu_1)$ and $a_1,a_2\geq 0$.\\
\itemnumber {ii} $\left\|Uf\right\|_\lapimN=\left\|f\right\|_\lapdoN$, for all $f\in \lap(S_1,\mu_1)$.\\
The max--linear isometry $U$ is called \textit{max--linear isomorphism} if it is onto.
\end{Def}

Consider a max--linear sub--space ${\cal F}\subset
L_+^\alpha(S_1,\mu_1)$ and a max--linear isometry $U:{\cal F} \to
L_+^\alpha(S_2,\mu_2)$.  Our goal in this section is somewhat
technical.  Namely, to characterize $U$ and also identify the largest max--linear 
sub--space ${\cal G} \subset L_+^\alpha(S_1,\mu_1)$, such that
$\filF\subset\calG$ and $U$ extends to $\calG$ uniquely as a
max--linear isometry.  This is done in Theorem \ref{hardin81thm:4.2p} below.
The proofs for all results in this section are given in Appendix~\ref{sec:proofMaxLinear}.

It is known that all linear isometries on $\la$ spaces for $\alpha\neq 2$ are related to a
 \textit{regular set isomorphism} (see~\cite{lamperti58isometries}). Regular set isomorphisms also play an important in the 
study of max--linear isometries. 
\begin{Def}\label{def:regular}
Let $(S_1,\calS_1,\mu_1)$ and $(S_2,\calS_2,\mu_2)$ be two measure spaces.
A set--mapping $T:\calS_1\to\calS_2$ is said to be a \textit{regular set isomorphism} if:\\
\itemnumber i For all $A\in\calS_1$, $T(S_1\backslash A) = T(S_1)\backslash T(A)\mod\mu_2$;\\
\itemnumber {ii} For disjoint $A_n$'s in $\calS_1$, $T(\cup_{n=1}^{+\infty}A_n) = \cup_{n=1}^{+\infty}T(A_n)\mod\mu_2$; \\
\itemnumber {iii} $\mu_2(T(A)) = 0$ if and only if $\mu_1(A) = 0$.
\end{Def}
\begin{Rem}
Regular set isomorphisms are mappings defined modulo null sets.  In the sequel, we often identify measurable sets
that are equal modulo null sets.
\end{Rem}
\noi The next properties follow immediately from the above definition: 

\itemnumber {iv} If $A_1,A_2\in\calS_1$ and $\mu_1(A_1\cap A_2) = 0$, then $\mu_2(T(A_1)\cap T(A_2)) = 0$.\\
\itemnumber {v} For all, not necessarily disjoint, $A_n \in \calS_1$, $n\in {\mathbb N}$, we have:
$$
T(\cup_{n=1}^\infty A_n) = \cup_{n=1}^\infty T(A_n)\ \ \ \mbox{ and }\ \ \ 
T(\cap_{n=1}^\infty A_n) = \cap_{n=1}^\infty T(A_n).$$

Any regular set isomorphism $T$ induces a canonical function mapping $Tf$, defined for all
measurable functions $f$, and such that $\{Tf \in B\} = T\{f \in B\}$, mod $\mu_2$, 
for all Borel sets $B\in \calB_{\mathbb R}$.  The resulting mapping is linear and also 
max--linear.  If $T$ is, in addition, measure preserving, then the induced mapping becomes 
a max--linear isometry. For more details, see Lemma~\ref{lem:maxLinearExtension} in Appendix~\ref{sec:proofMaxLinear} or Doob~\cite{doob53stochastic}.  
The next result shows that any max--linear isometry, which maps the identity 
function $\ind$ to the identity function $\ind$, is induced by a measure preserving regular set isomorphism.

\begin{Thm}\label{hardin81thm:2.2p}
Suppose $\alpha>0$. Let $\filF$ be a max--linear sub--space of $\lapdo$ and $U:\filF\to\lapim$ be a max--linear isometry. If $\inddo\in\filF$ and $U\inddo = \indim$, then $Uf= Tf$ for all $f\in\filF$, where:

\itemnumber i $T$ is induced by a measure preserving regular set isomorphism from $\sigma(\filF)$ onto $\sigma(U(\filF))$,\\
\itemnumber {ii} $T$ is a max--linear isometry from $\lapdoS$ onto $\lapimS$, and\\
\itemnumber {iii} $T$ is the unique extension of $U$ to a max--linear isometry from $\lapdoS$ to $\lapim$.
\end{Thm}

Not all max--linear isometries are directly induced by regular set isomorphisms. We will
show next, however, that every max--linear isometry can be related to a regular set 
isomorphism.

\begin{Def} Let $F$ be a collection of functions in $\lap(S,\mu)$.

\itemnumber i The \textit{ratio $\sigma$--field} of $F$, written $\rho(F)\defe\sigma\left(\left\{f_1/f_2, f_1,f_2\in F\right\}\right)$, is defined as the $\sigma$--field generated by ratio of functions in $F$, where the ratios take values in the extended interval $[0,\infty]$;\\
\itemnumber {ii}  The \textit{positive ratio space} of $F$, written $\ratios F$, is defined as $\lap(S,\rho(F),\mu)$.
\\
\itemnumber {iii} The \textit{extended positive ratio space} of $F$, written $\eratios F$, is defined as the class of all functions in $\lap(S,\mu)$ that have the form $rf$, where $r$ is non-negative $\rho(F)$-measurable and $f\in F$. 
\end{Def}

\noindent In the following lemma, we present some important properties of the ratio $\sigma$--fields.
\begin{Lem}\label{lem:ratio1} 
For any non--empty class of functions $F \subset \laps$, we have
 $\rho(F) = \rho(\cmsspan(F))\subset\sigma(F)$. If, in addition, $\ind_S\in F$, then $\rho(F) = \sigma(F)$.
\end{Lem}
 
Before introducing the main result of this section, we need some auxiliary results about the notion of \textit{full support}.

\begin{Def}\label{def:fullSupport}
Let $(S,\mu)$ be a measurable space and $F$ be a collection of
measurable real-valued functions on $(S,\mu)$. A measurable function
$f_0$ is said to have \textit{full support} w.r.t. $F$ if
$\mu(\supp(g)\setminus \supp(f_0)) = 0$ for all $g\in F$, where
$\supp(f)\defe \{f\neq 0\}$. If, in addition, $f_0\in F$, we then
write $\supp(F) = \supp(f_0)$.
\end{Def}
\begin{Rem}
Note that the definition of full support is modulo $\mu$-null sets and the definition of $\supp(F)$ is independent of the choice of $f_0\in F$. Also, our definition of $\supp(F)$ requires implicitly that $F$ contains a function $f_0$ of full support.
\end{Rem}

\begin{Lem}\label{lem:fullSupport1}
Let $\filF$ be a max--linear sub--space of $\laps$. If $\filF$ is separable or $\mu$ is $\sigma$-finite, then there exists a function of full support in $\filF$.
\end{Lem}
\begin{Lem}\label{lem:fullSupport2}
Let $\filF$ be a max--linear sub--space of $\lapdo$ and let $U:\filF\to\lapim$ be a max--linear isometry. Assume that the measures $\mu_1$ and $\mu_2$ 
are $\sigma$--finite. If $f_0$ has full support in $\filF$, then $Uf_0$ has full support in $U(\filF)$.
\end{Lem}

\noi We now present the main result of this section.
\begin{Thm}\label{hardin81thm:4.2p}
Suppose $\alpha>0$ and let $\calF$ be a max--linear sub--space of $\lap(S_1,\mu_1)$.  Suppose also that
$\supp({\calF}) = S_1$. If $\mu_1$ is $\sigma$-finite and $U:\filF\to\lapim$ is a max--linear isometry, then:

\noi\itemnumber i
$U$ has a unique extension to a max--linear isometry $\widebar U$, defined on 
$\eratiosf$ to $\lap(S_2,\mu_2)$.  Moreover, $\widebar U$ is also onto $\eratiosuf \subset\lap(S_2,\mu_2)$ and
\equh\label{hardin81thm:4.2peq:1}
\widebar U(rf) = (Tr)(Uf),\ \ \ \mbox{ for all } r\in\ratiosf\,,f\in\filF\,,
\eque
where the function mapping $T:\ratiosf\to\ratiosuf$ is induced by a regular set isomorphism of $\rho(\filF)$ onto $\rho(U(\filF))$. \\
\itemnumber {ii} For all $f\in \calF,$ we have 
\equh\label{hardin81thm:4.2peq:2}
(Uf)^\alpha d\mu_2 = d\mu_{1,f}\circ{T\inv}\,,
\eque
where $d\mu_{1,f} = f^\alpha d\mu_1$.
\end{Thm}
\begin{Rem} Equality \eqref{hardin81thm:4.2peq:2} means that the two measures are identical on
the $\sigma$--field $\rho(U(F))$, i.e.\ $\int_A (Uf)^\alpha d\mu_2 = \mu_{1,f}\circ T^{-1}(A)$,
for all $A \in \rho(U(F))$.  In the sequel, we will interpret equalities between measures defined on 
different $\sigma$--fields as equality of their corresponding restrictions to the largest common
$\sigma$--field.  Note that in general $(Uf)^\alpha$ in \eqref{hardin81thm:4.2peq:2} does not 
necessarily equal the Radon--Nikodym derivative $d(\mu_{1,f}\circ T^{-1})/d\mu_2$ since the
$\sigma$--field $\rho(U(F))$ is typically rougher than ${\cal B}_{S_2}$. This is why $U$ may not have
a unique extension to $L_+^\alpha(S_2,\mu_2)$, in general.  See Remark 3.2(c) in 
Rosi\'nski~\cite{rosinski06minimal} for a detailed discussion.
\end{Rem}
Recall the notion of equivalence in measure of two $\sigma$--fields, defined on the same measure
space $(S,{\cal S},\mu)$.  Namely, for two $\sigma$--fields $\calA,\ {\calB} \subset {\cal S}$, we
write $\calA\sim\calB\mod\mu$, if for any $A \in \calA$ ($B\in\calB$, respectively), there exists 
$B \in \calB$ ($A\in\calA$, respectively) such that $\mu(A\Delta B) = 0$. 
The following result will be used in the next section.

\begin{Lem}\label{lem:ratio2}
Let $F$ be a class of functions in $\laps$. Suppose there exists $f_0\in F$ with full support in $F$.  
If $S = \supp(f_0) \equiv \supp(F)$ and
if $\rho(F) \sim \calB_S \mod\mu$, then $\eratios F = \laps$.
\end{Lem}

This result and Theorem \ref{hardin81thm:4.2p}, provide sufficient conditions for a max--linear isometry $U$, defined on
$F$, to extend uniquely to the entire space $\laps$. 

\section{Minimal Representations for $\alpha$--Fr\'echet Processes}\label{sec:minimal}

Let $\{f_t^{(i)}\}_{t\in T} \subset \lap(S_i,\mu_i),\ i=1,2$ be two spectral representations for the $\alpha$--Fr\'echet
process $X = \{X_t\}_{t\in T}$. Recall that for all $t_j \in \bbR,\ c_j \ge 0,\ 1\le j\le n,$ we have
$$
\P\{X_{t_j}\le c_j^{-1},\ 1\le j\le n\} = \int_{S_1} {\Big(} \bigvee_{j=1}^n c_j f_{t_j}^{(1)} {\Big)}^\alpha d\mu_1
= \int_{S_2} {\Big(} \bigvee_{j=1}^n c_j f_{t_j}^{(2)} {\Big)}^\alpha d\mu_2.
$$
One can thus define the following natural max--linear isometry: 
\equh\label{eq:Ucanonical}
 U: \cmsspan\{f_t^{(1)}\}_{t\in T}
 \to \cmsspan\{f_t^{(2)}\}_{t\in T}\,,\ \mbox{ with }\ U f_t^{(1)} := f_t^{(2)},\ \mbox{ for all }t\in T.
\eque 
In the sequel, $U$ will be called the {\it relating max--linear isometry} of the two representations.   
Our goal in this section is to provide convenient representations for the max--linear isometry $U$.

For any standard Lebesgue space $(S,\mu)$, we have that $\{f_t\}_{t\in T}\subset L_+^\alpha(S,\mu)$ is separable, and hence by
Lemma \ref{lem:fullSupport1}, the max--linear space ${\cal F}=\cmsspan(f_t,\ t\in T)$ contains a function with 
{\it full support}.  Therefore, by convention, we define the support of $\{f_t\}_{t\in T}$ as follows:
$$
{\rm supp}\{ f_t,\ t\in T\} := \supp({\cal F}) \equiv {\rm supp} {\Big(} \cmsspan(f_t,\ t\in T) {\Big)}.
$$
In view of Theorem \ref{hardin81thm:4.2p}, one can readily represent the max--linear isometry $U$ 
in \eqref{eq:Ucanonical} in terms of a regular set isomorphism.  The latter mapping however is a set--mapping rather than point mapping.
It is desirable to be able to express $U$ via measurable point mappings. Unfortunately, in general such point mappings may not be unique.
In order to have a unique point mapping relating the two representations, we need to impose further \textit{minimality condition} 
on the spectral representations. The following definition is as in Rosi\'nski~\cite{rosinski95structure} (see also~\cite{hardin82spectral}).

\begin{Def}\label{def:minimality} A spectral representation
 $\indt f\subset\lap(S,\mu)$ of an $\alpha$--Fr\'echet process is said to be {\it minimal} if:\\
\noindent \itemnumber i $\supp\{f_t: t\in T\} = S \quad \mu\ae$, and\\
\itemnumber {ii} for any $B\in\calB_S$, there exists $A\in\rho(\{f_t: t\in T\})$ such that $\mu(A\Delta B) = 0$.
\end{Def}
We shall also consider minimal representations with \textit{standardized support} defined as follows. 
\begin{Def}\label{def:standardizedSupport}
A minimal representation $\indt f\subset \laps$ has {\it standardized support} 
if, up to $\mu$-null sets:\\
\noindent
\itemnumber {i} $S\subset(0,1)\cup\mathbb N$,\\
\itemnumber {ii} $S\cap(0,1) = \emptyset$ or $(0,1)$ and $\mu|_{(0,1)}$ is the Lebesgue measure, \\
\itemnumber{iii} $S\cap\mathbb N = \emptyset$, $\mathbb N$ or $\{1,\cdots,N\}$, where $N\in\mathbb N$ and $\mu|_{S\cap\mathbb N}$ is the counting measure.\\
\noindent Let $(\ssS,\ssL)$ denote the standard support with $I=0$ or $1$ respectively according to the two cases in (i) 
and $N=0,N = \infty$ or $N\in\mathbb N$ respectively according to the three cases in (ii), e.g. $S_{0,\infty} = \mathbb N$ and $S_{1,N} = (0,1)\cup\{1,\dots,N\}$.
\end{Def}

We now show that any spectral representation of an $\alpha$--Fr\'echet process can be transformed into a minimal one with standardized support.
\begin{Thm}\label{thm:standardized}
Every $\alpha$--Fr\'echet process satisfying Condition S has a minimal representation $\indt f$ with standardized support $(\ssS,\ssL)$. That is
\equh\label{rep:minRepSS}
\indtauto X \eqd {\Big\{}\Eintt_\ssS f_t(s)M_\alpha(ds) {\Big\}}_{t\in T}\,,
\eque
where $M_\alpha$ is the $\alpha$--Fr\'echet random sup--measure with control measure $\ssL$.
\end{Thm}

\begin{proof} By Proposition~\ref{prop:conditionS}, one can 
let $G = \{g_t\}_{t\in T}\subset \lap((0,1),\calB_{(0,1)},ds)$ be a spectral representation
of the process in question, where $ds$ is the Lebesgue measure on $(0,1)$. First, we study the 
ratio $\sigma$--field generated by $G$. Let $\calG = \cmsspan\{g_t,t\in T\}$ and, in view of 
Lemma~\ref{lem:fullSupport1}, let $g\in\calG$ have full support in $\calG$. By Lemma~\ref{lem:ratio1}, we have
$\rho(G) = \rho(\calG)$. Without loss of generality we assume $\supp(g) = \supp(\calG) = (0,1)$ and 
$\left\|g\right\|_\alpha = 1$. Define a new measure $\mu$ on the space $((0,1),\rho(\calG))$ by 
setting $d\mu(s) = g(s)^\alpha ds$. Since $\mu$ is a probability measure, the measure space 
$((0,1),\rho(\calG),\mu)$ has at most countably many (equivalence classes of) atoms. With some 
abuse of notation, we represent them as $A_1,A_2,\dots,A_N$, where $N=0$ means no atoms, $N \in 
\mathbb N$ for finite number of atoms, and $N = \infty$ when countably infinite number of atoms are 
present. Set $A=\cup_{n=1}^N A_n$ and $a_i = \mu(A_i)\,,1\leq i\leq N$. 

Next, we define a regular set isomorphism $T_r$ of measure space $((0,1),\rho(\calG),\mu)$ onto 
measure space $(\ssS,\ssB,\ssL)$ considered in Definition~\ref{def:standardizedSupport}. For the 
atoms, define $T_r^N(A_n) = \{n\}, n\leq N,n\in\mathbb N$. For the non--atomic subset $A_0\equiv 
(0,1)\setminus A$, let $\calS_0 = \rho(\calG)\cap A_0 =  \{B\cap A_0, B\in \rho(\calG)\}$ and let 
$\mu_i$ be the restriction of $\mu$ to $A_i, i=0,\dots,N$. 

The case $a_0 = 0$ is trivial since then $\mu(A_0) = \mu_0(A_0) = 0$ and we can simply ignore 
$(A_0,\calS_0,\mu_0)$. 
We thus suppose that $a_0>0$ and observe that $(A_0,\calS_0,\mu_0)$ is an non--atomic separable 
measurable space (see p167 in~\cite{halmos50measure}) with total mass $\mu(A_0) = 1-\sum_{n=1}^Na_n 
\equiv a_0$. Indeed, the separability of $(A_0,\calS_0,\mu_0)$ is due to the fact that $\calG$ restricted on 
$A_0$ is separable. 

Now, Theorem 41.C in Halmos \cite{halmos50measure} implies that there is a measure preserving regular 
set isomorphism, i.e., a {\it measure algebra isomorphism} $T_r^I$ from $(A_0,\calS_0,\mu_0)$ {\it 
onto} $((0,1),\calB_{(0,1)},a_0ds)$. By combining the definitions of $T_r^N$ on all atoms $A_i,\ 1\le i \le N$ and $T_r^I$ on 
$(A_0,\calS_0,\mu_0)$, we thus obtain a regular set isomorphism $T_r\defe T^I_r+T^N_r$ from 
$((0,1),\rho(\calG),\mu)$ {\it onto} $(\ssS,\ssB,\ssL)$. Note that $T_r$ is not necessarily measure 
preserving.

By using $T_r$, we construct next the desired minimal representation with standardized support. 
Define
\equh\label{hardin82thm:1.1peq:isometry}
f_t(s) = T_r(g_t/g)(s)\left(a_0^{1/\alpha}\ind_{(0,1)}(s)+\sum_{n=1}^Na_n^{1/\alpha}\ind_{\{n\}}(s)\right)\,,
\eque
where $T_r$ is the canonical map on measurable functions induced by the constructed isomorphism 
(see Lemma~\ref{lem:maxLinearExtension} or p452-454 \cite{doob53stochastic}) from $\lap((0,1),
\rho(\calG),\mu)$ onto $\lap(\ssS,\ssL)$. 
We claim that $\fft$ is a minimal representation with standardized support. It is clearly a 
spectral representation, since, for any $m\in\mathbb N, t_i\in T,ci_i>0,1\leq i\leq m$, 
\eqnhn
{\Big\|}\bigvee_{i=1}^mc_if_{t_i}{\Big\|}_{\lap(\ssS,\ssL)}^\alpha
& = & 
{\Big\|}\bigvee_{i=1}^mc_iT_r(g_{t_i}/g)
{\Big(}a_0^{1/\alpha}\ind_{(0,1)}+\sum_{n=1}^Na_n^{1/\alpha}\ind_{\{n\}}{\Big)}
{\Big\|}_{\lap(\ssS,\ssL)}^\alpha \nonumber\\
& = & {\Big\|}a_0^{1/\alpha}\bigvee_{i=1}^mc_iT_r(g_{t_i}/g){\Big\|}_{\lap(0,1)}^\alpha + \sum_{n=1}^N{\Big|}a_n^{1/\alpha}\bigvee_{i=1}^mc_iT_r(g_{t_i}/g)(n){\Big|}^\alpha \nonumber\\
& = & {\Big\|}\bigvee_{i=1}^mc_ig_{t_i}/g{\Big\|}_{L^\alpha(A_0,\mu_0)}^\alpha + \sum_{n=1}^N{\Big\|}\bigvee_{i=1}^mc_ig_{t_i}/g{\Big\|}_{L^\alpha(A_n,\mu_n)}^\alpha \label{eq:scaling}\\
& = & {\Big\|}\bigvee_{i=1}^mc_ig_{t_i}/g{\Big\|}_{L^\alpha((0,1),\mu)}^\alpha = {\Big\|}\bigvee_{i=1}^mc_ig_{t_i}{\Big\|}_{\lap(0,1)}^\alpha \nonumber
\eqnen
where \eqref{eq:scaling} follows from the fact that $T_r^I$ is a measure preserving regular set isomorphism of $A_0$ 
onto $(0,1)$ and since $T_r^N$ maps atoms to integer points in a one-to-one and onto manner. Indeed, restricted on each
$A_i, 0\le i\le N$, $a_i^{1/\alpha}T_r$ is a max--linear isometry satisfying
\eqnh
\Big\|a_i^{1/\alpha}T_r\ind_{A_i}\Big\|_{L^\alpha(T_r(A_i),\ssL)}^\alpha & = & \Big\|a_i^{1/\alpha}\ind_{T_r{A_i}}\Big\|_{L^\alpha(T_r(A_i),\ssL)}^\alpha\\
& = & a_i\ssL(T_rA_i) = \mu(A_i) = \left\|\ind_{A_i}\right\|_{L^\alpha(A_i,\mu_i)}^\alpha\,.
\eqne
We will complete the proof by verifying the minimality of $\{f_t\}_{t\in T}$ (by Definition~\ref{def:minimality}). 
Let $\filF$ denote $\cmsspan\{f_t,t\in T\}$ and note that $g\in \calG = 
\cmsspan\{g_t,\ t\in T\}$. Since $T_r(g/g) = \ind_{S_{I,N}}$, 
by~\eqref{hardin82thm:1.1peq:isometry}, we obtain that
\equh\label{eq:ssf}
\ssF(s) \defe 
a_0^{1/\alpha}\ind_{(0,1)}(s)+\sum_{n=1}^Na_n^{1/\alpha}\ind_{\{n\}}(s) \ \ \mbox{ belongs to }\ \  \filF.
\eque
This implies $\supp(\ssF) = \supp(\filF) = \ssS$, and whence (i) in Definition~\ref{def:minimality} holds. 
To verify (ii), observe that by \eqref{hardin82thm:1.1peq:isometry} and Lemma~\ref{lem:maxLinearExtension}, 
$f_1/f_2 = T_r(g_1/g)/T_r(g_2/g) = T_r(g_1/g_2)$ for all $g_1,g_2\in\calG$.  Therefore $T_r(\rho(\calG)) \equiv \rho(\calF)$, and since, as shown above, the regular set isomorphism $T_r$ maps $\rho(\calG)$ {\it onto} $\calB_{\ssS}$, it follows that
(ii) holds. \ifthenelse{\boolean{qedTrue}}{\qed}{}
\end{proof}
\begin{Rem} Theorem \ref{thm:standardized} shows the existence of minimal representations with standardized support. One
 can have many minimal representations whose supports are not necessarily standardized in the same way.  For example, in 
 the proof of Theorem \ref{thm:standardized}, we could define
$\tilde\lambda_{I,N}$ on $S_{I,N}$ so that restricted on the atoms $A_i$, $1\le i \le N$, we have
$d\tilde\lambda_{I,N} = a_i^{1/\alpha}d\ssL$.  In this case, one obtains a finite measure $\tilde\lambda_{I,N}$ on $S_{I,N}$ as discussed in Rosi\'nski \cite{rosinski94uniqueness} (p. 626) for the case of symmetric $\alpha$--stable processes.
Our measure $\lambda_{I,N}$ may be infinite, since it is a counting measure on the atoms.

\end{Rem}
\begin{Rem}
Theorem~\ref{thm:standardized} can be seen as a generalization of Theorem 4.1 in de Haan and Pickands III~\cite{dehaan86stationary}. Instead of minimal representation, \textit{proper representation} is involved therein. A spectral representation is proper if the spectral functions $\indt f$ satisfy (i) $\supp\{f_t\,,t\in T\} = S\,,\mu\ae$ and (ii) $\forall B\in\calB_S$, either there exists $A\in \rho(\{f_t\,,t\in T\})$ such that $\mu(A\Delta B) = 0$ or there exists an atom $A\in\rho(\{f_t\,,t\in T\})$ such that $\mu(B\cap A)>0$.
 This definition is closely related to our definition of minimality, in the sense that any proper representation can be transformed into a minimal one. Indeed, this essentially involves contracting the atoms to points as in the proof of Theorem~\ref{thm:standardized}.
\end{Rem}

Consider the \textit{canonical} max--linear isometry $U$ relating two
spectral representations as in~\eqref{eq:Ucanonical}.
Theorem \ref{hardin81thm:4.2p} implies that $U$ extends uniquely to a max--linear isometry 
$U:\eratios{\filF\topp1}\to\eratios{\filF\topp2}$ between extended positive ratio spaces, where
$\filF\topp i = \cmsspan\{f_t\topp i:t\in T\}\,,i=1,2$. Now, if the first spectral representation
$\indt{f\topp1}$ is \textit{minimal}, then by Lemma~\ref{lem:ratio2},
$\eratios{\filF\topp1} = \lap(\ssS,\ssL)$. In this case, one can also
represent $U$ in terms of \textit{measurable point mappings}. This
\textit{point mapping representation} is developed in the following
result. It will be essential for our studies in Sections
\ref{sec:classification} and \ref{sec:stationary}.

\begin{Thm}\label{thm:relation}
Let $\indt f\subset\lap(\ssS,\ssL)$ and $\indt g\subset\laps$ be two
spectral representations of an $\alpha$--Fr\'echet process $\indt
X$. Let $U$ be the relating max--linear isometry of $\indt f$ and
$\indt g$. If $\indt f$ is minimal and $\indt g$ is arbitrary,
then\\ \itemnumber i $U$ can be uniquely extended to
$\lap(\ssS,\ssL)$; \\ \itemnumber {ii} $U$ can be represented by
measurable functions $\Phi:S\to\ssS$ and $h:S\to\mathbb
R_+\setminus\{0\}$, such that $\Phi$ is onto, and the following
statements hold: 
\equh\label{eq:pointRep1} g_t(s) = Uf_t(s) =
h(s)\left(f_t\circ\Phi\right)(s)\,, \quad\mu\ae \,,
\eque 
and
\equh\label{eq:pointRep2} d\ssL = d\left(\mu_h\circ\Phi\inv\right),
\eque 
where $d\mu_h(s) = h(s)^\alpha d\mu$. $\Phi$ is unique modulo
$\mu$.
\end{Thm}
\begin{proof}
Let $F$ and $G$ denote $\indt {f}$ and $\indt g$ respectively.  By Theorem~\ref{hardin81thm:4.2p},  
there exists a regular set isomorphism $T_r$ from $\ssB$ {\it onto} $\rho(G)$ such that
\[
g_t(s) = Uf_t(s) = (T_rf_t)(s) \left(\frac{Uf_0}{T_r f_0}\right)(s),\ \ \ \mu\ae\,,\forall t\in T\,,
\]
for some function with full support $f_0 \in \cmsspan\{f_t,\ t\in T\}$. In the last relation we used the facts
that $T_r(1/f_0)= 1/T_r(f_0)$ and $T_r(f_t/f_0) = T_r(f_t)/T_r(f_0)$ (Lemma \ref{lem:maxLinearExtension}).  
Moreover, we have that
\equh\label{eq:uf0a}
\left(Uf_0\right)^\alpha d\mu = d(\mu_{1,f_0}\circ{T_r\inv}) = \left(T_r f_0 \right)^\alpha d\left(\ssL\circ T_r\inv\right)  \,,\mu\ae\,.
\eque
By Theorem 32.5 in Sikorski~\cite{sikorski64boolean}, the regular set isomorphism $T_r$ can be induced by a point mapping $\Phi$ from $S$ onto $\ssS$ such that 
$T_rf = f\circ \Phi$, for all measurable functions $f$ defined on $S_{I,N}$. Moreover, $\Phi$ is unique modulo $\mu$. Note that in general
$\Phi$ is not one-to-one, because of the possible presence of atoms in $(S,\rho({\cal G}),\mu)$.  To show that \eqref{eq:pointRep2} is true, let
\[
  \widetilde h(s) = \frac{Uf_0}{T_r f_0}(s) = \frac{Uf_0}{f_0\circ \Phi}(s)\,.
\]
Note that by Lemma~\ref{lem:fullSupport2}, $\widetilde h(s)>0\,,\mu\ae$. Put 
\equh\label{eq:hs}
h(s) = \left\{
\begin{array}{l@{\mbox{ if }}l} \widetilde h(s) & \widetilde h(s)>0\\ 1 & \widetilde h(s) = 0\end{array}
\right.
\mand d\mu_h = h^\alpha d\mu.
\eque
Observe that $h$ is a measurable function from $S$ to $\mathbb R_+\setminus\{0\}$. Thus, relation~\eqref{eq:pointRep2} follows by~\eqref{eq:hs} and~\eqref{eq:uf0a}.
This completes the proof.\ifthenelse{\boolean{qedTrue}}{\qed}{}
\end{proof}

\begin{Rem}\label{rem:Phi-non-sing}
Relation \eqref{eq:pointRep2} and the fact that $h(s)>0$ for all $s$ imply that $\mu\circ 
\Phi^{-1}\sim \lambda_{I,N}$. 
\end{Rem}

\noindent 
Now, if both representations in Theorem~\ref{thm:relation} are minimal, we have the following:
\begin{Coro}\label{coro:uniquePointMapping}
If $\indt {f\topp i}\,,i=1,2$ are two minimal representations of an $\alpha$--Fr\'echet process $\indt X$ with standardized support $(S_{I_i,N_i},\lambda_{I_i,N_i})\,,i=1,2$, then the relating max--linear isometry $U$ from $\lap(S_{I_1,N_1},\lambda_{I_1,N_1})$ onto $\lap(S_{I_2,N_2},\lambda_{I_2,N_2})$ is determined by, unique modulo $\lambda_{I_2,N_2}$, functions $\Phi: S_{I_2,N_2}\to S_{I_1,N_1}$ and $h: S_{I_2,N_2}\to\mathbb R_+\setminus\{0\}$ such that $\Phi$ is one-to-one and onto and, for each $t\in T$,
\equh\label{eq:minimalRelation1}
f_t\topp2(s) = Uf_t\topp1(s) = h(s)\left(f_t\topp1\circ\Phi\right)(s)\,, \quad \lambda_{I_2,N_2}\ae
\eque
and
\equh\label{eq:minimalRelation2}
\frac{d(\lambda_{I_1,N_1}\circ\Phi)}{d\lambda_{I_2,N_2}}(s) = h(s)^\alpha,\quad \lambda_{I_2,N_2}\ae\,.
\eque
\end{Coro}

\noindent
An important consequence of Corollary~\ref{coro:uniquePointMapping} is the following.
\begin{Coro}\label{coro:uniqueness}
Let $\indt {f\topp i}\,,i=1,2$ be as in Corollary~\ref{coro:uniquePointMapping}. Then 
\[
I_1 = I_2 = I\quad\mbox{and}\quad N_1 = N_2 = N\,.
\]
Moreover, the relating max--linear isometry $U:\lap(\ssS,\ssL)\to\lap(\ssS,\ssL)$ satisfies\\
\itemnumber i if $I = 1$, then $\forall f\in\lap(0,1)$,
\equh\label{eq:ufi}
Uf = \left(\dfrac\lambda{\Phi_I}\lambda\right)^\alpha(f\circ\Phi_I)\,,\lambda\ae\,,
\eque
where $\lambda$ is the Lebesgue measure on $(0,1)$, $\Phi_I$ is a point map from $(0,1)$ onto $(0,1)$, and\\
\itemnumber {ii} if $N \neq 0$, then $\forall f\in\lap(\ssS\cap\mathbb N,\ssL)$,
\equh\label{eq:ufn}
Uf = f\circ\Phi_N\,,
\eque
where $\Phi_N$ is an automorphism of $\ssS\cap\mathbb N$.
\end{Coro}
\begin{proof} 
We start by recalling that $U$ is induced by $T_r$, which is an one-to-one isomorphism modulo $\ssL$-null sets from
$\calB_{S_{I_1,N_1}}$ onto $\calB_{S_{I_2,N_2}}$ (by Theorem~\ref{hardin81thm:4.2p}).
Since $T_r$ is a regular set isomorphism, one has that for all $A,B\in\calB_{S_{I_1,N_1}}$,
\[
\lambda_{I_1,N_1}(A)\lambda_{I_1,N_1}(B\setminus A) = 0\Leftrightarrow \lambda_{I_2,N_2}(T_rA)\lambda_{I_2,N_2}(T_rB\setminus T_rA) = 0\,.
\]
Thus $T_r$ maps \textit{atoms} to \textit{atoms} and non-atomic sets to non-atomic sets. Hence,
\[
T_r\left(\calB_{S_{I_1,N_1}}\cap(0,1)\right) \subset \calB_{S_{I_2,N_2}}\cap(0,1) \mbox{ and }
T_r\left(\calB_{S_{I_1,N_1}}\cap\mathbb N\right) \subset \calB_{S_{I_2,N_2}}\cap\mathbb N\,.
\]
Since $T_r$ is onto, we also have that
\[
T_r\left(\calB_{S_{I_1,N_1}}\cap(0,1)\right) = \calB_{S_{I_2,N_2}}\cap(0,1) \mbox{ and } T_r\left(\calB_{S_{I_1,N_1}}\cap\mathbb N\right) = \calB_{S_{I_2,N_2}}\cap\mathbb N\,.
\]
This implies that $I_1=I_2$. Moreover, since $T_r$ is one--to--one and onto, we have $N_1 = N_2$. This also shows that $T_r:\ssS\cap\mathbb N\to\ssS\cap\mathbb N$ is a bijection where $I \defe I_1=I_2$ and $N \defe N_1=N_2$. By Corollary~\ref{coro:uniquePointMapping}, it follows that (i) and (ii) holds. Note that in (ii) we have simpler formula for $Uf$. This is because that on the discrete part $\ssS\cap\mathbb N$, the function $h(s)$ defined in~\eqref{eq:minimalRelation2} equals 1.\ifthenelse{\boolean{qedTrue}}{\qed}{}
\end{proof}
\begin{Rem}
Theorem~\ref{thm:relation} and Corollary~\ref{coro:uniquePointMapping} are valid even if the minimal representations therein do not have standardized support (see Theorem~4.1 and Theorem~4.2 in~\cite{dehaan86stationary} for results on discrete processes; see also Theorem~2.1 in Rosi\'nski~\cite{rosinski95structure} for analogous result in the sum--stable setting). The advantage of having minimal representation with \textit{standardized support} is shown in Corollary~\ref{coro:uniqueness} and further exploited in the next section.
\end{Rem}
\section{Classification of $\alpha$--Fr\'echet Processes}\label{sec:classification}

We now apply the abstract results on max--linear isometries and minimal representations
to classify $\alpha$--Fr\'echet processes.  The first classification result is an immediate
consequence of the notion of minimal representation with standardized support and it applies to 
general max--stable processes.

\subsection{Continuous--discrete decomposition}\label{sec:continuousDiscreteDecomposition}

Consider an $\alpha$--Fr\'echet process $X=\{X_t\}_{t\in T}$, which has a minimal representation with standardized support $\{f_t\}_{t\in T}\subset L_+^\alpha(S_{I,N},\lambda_{I,N})$.
By Corollary~\ref{coro:uniqueness}, the support $(\ssS,\ssL)$ is unique.  We therefore call
$S_{I,N}$ the {\it standardized support} of $X$ and focus on the {\it continuous} and {\it discrete} parts of $S_{I,N}$, respectively:
\[
 S_I \defe \ssS\cap(0,1),\quad\mbox{ and }\quad S_N\defe \ssS\cap\mathbb N.
\]
Let $f\tpd I_t = f_t\ind_{S_I},$ and $f\tpd N_t = f_t\ind_{S_N}$ be the restrictions of
the $f_t$'s to $S_I$ and $S_N$, respectively.
One can write:
\equh\label{eq:continuousDiscreteDecomposition}
 \indt X \defd \left\{X^I_t \vee X^N_t\right\}_{t\in T},
\eque
where
\equh\label{eq:continuousDiscreteDecomposition2}
X_t^I := \Eintt_{S_I}f^I_t(s)M_\alpha(ds)\quad \mbox{ and }\quad
X^N_t := \Eintt_{S_N}f^N_t(s)M_\alpha(ds)\,,
\eque
are two independent $\alpha$--Fr\'echet processes.  
The following result shows that the decomposition \eqref{eq:continuousDiscreteDecomposition} does
not depend on the choice of the representation $\{f_t\}_{t\in T}$. 

\begin{Thm}\label{thm:continuousDiscreteDecomposition} Let $\indt X$ be an $\alpha$--Fr\'echet process
with minimal representation of standardized support $\{f_t\}_{t\in T}\subset L_+^\alpha(S_{I,N},\lambda_{I,N})$.
Then:\\
\itemnumber {i} The decomposition~\eqref{eq:continuousDiscreteDecomposition} is unique in distribution.

\itemnumber {ii} The processes $X^I=\indt {X^I}$ and $X^N=\indt {X^N}$ are independent and they have
standardized supports $S_I$ and $S_N$, respectively.

\itemnumber {iii} The functions $\{f_t^I\}_{t\in T} \subset L_+^\alpha(S_I,\lambda_I)$ and
$\{f_t^N\}_{t\in T} \subset L_+^\alpha(S_N,\lambda_N)$ provide minimal representations for the processes $X^I$ and
$X^N$, respectively.
\end{Thm}
\begin{proof}
To prove {\it (i)}, suppose $\indt g\subset\lap(\ssS,\ssL)$ is another minimal representation of $X$ with standardized support and consider the decomposition $\indt X\eqd \left\{Y^I_t\vee Y^N_t\right\}_{t\in T}$, where  
\[
Y_t^I := \Eintt_{S_I}g^I _t(s)M_\alpha(ds)\quad \mbox{ and }\quad
Y^N_t := \Eintt_{S_N}g^N _t(s)M_\alpha(ds)\,,\forall t\in T\,.
\]
By Corollary~\ref{coro:uniqueness}, the relating max--linear isometry
$U$ of $\indt f$ and $\indt g$ is such that for all $t\in T$,
$U(f_t^I) = g_t^I$ and $U(f_t^N) = g_t^N$. Moreover, $U$ remains a
max--linear isometry when restricted to $S_I$ and $S_N$, and hence
\[
\indt {X^I} \eqd \indt {Y^I}\quad \mbox{ and }\quad\indt {X^N} \eqd \indt {Y^N}\,.
\]
The last two relations imply that the
decomposition~\eqref{eq:continuousDiscreteDecomposition} does not depend on
the choice of the representation.  The components $\indt{\xI}$ and $\indt{\xN}$ are independent
since they are defined by extremal integrals over two disjoint sets $S_I$ and $S_N$.
The minimality of $\indt f$ implies the minimality of $\indt{f^I}$  and $\indt {f^N}$, restricted to
$S_I$ and $S_N$, respectively.  This completes the proof, since the supports $S_I$ and $S_N$ of
$\indt{f^I}$ and $\indt {f^N}$ are standardized (Definition \ref{def:standardizedSupport}). 
\end{proof}

The processes $\indt{\xI}$ and $\indt{\xN}$ in the Decomposition \eqref{eq:continuousDiscreteDecomposition} will be referred to as
the {\it spectrally continuous} and {\it spectrally discrete} components of $X$, respectively.  The next result clarifies
further their structure.

\begin{Coro}\label{coro:continuousDiscreteDecomposition} 
Let $\indt f$ and $\indt g$ be two minimal representations with standardized support of an
$\alpha$--Fr\'echet process $\indt X$.  Then, the relating max--linear
isometry $U$ of these representations, has the form
\equh\label{eq:ufti}
 Uf_t^ I = \left(\dfrac\lambda{\Phi_I}\lambda\right)^{1/\alpha}(f_t^
I\circ\Phi_I) = g_t^ I\quad\mbox{and}\quad Uf_t^ N = f_t^ N\circ\Phi_N
= g_t^ N\,, \lambda\ae, \ \ \forall t\in T,
\eque
where $\Phi_I$ is a point mapping from $S_I$ onto $S_I$ and $\Phi_N$ is a 
permutation of $S_N$ (a one-to-one mapping from $S_N$ onto $S_N$).
\end{Coro}

\noindent
The proof is an immediate consequence of Relations~\eqref{eq:ufi} and~\eqref{eq:ufn} above.  This result shows
that the discrete component of an $\alpha$--Fr\'echet process has an interesting invariance property.  
Namely, suppose that $X$ has a non--trivial discrete 
component $X^N = \{X_t^N\}_{t\in T}$. By Corollary \ref{coro:continuousDiscreteDecomposition}, there exists 
a {\it unique} set of functions $t\mapsto \phi_t(i),\ i\in S_N,\ t\in T$, such that:
{\it (i)} ${\rm supp}\{ \phi_t,\ t\in T\} \equiv S_N$, {\it (ii)} $\rho\{\phi_t,\ t\in T\}  = \calB_\ssS \equiv 2^{S_N}$ and
{\it (iii)} $\sum_{1\le i \le N} \phi_t(i)^\alpha <\infty$, for all $t\in T$ and
$$
 \{X_t^N\}_{t\in T} \eqd \bigvee_{i=1}^N \phi_t(i) Z_i,
$$
where $Z_i,\ 1\le i\le N$ are independent standard $\alpha$--Fr\'echet random variables.
The functions $t\mapsto \phi_t(i),\ 1\le i\le N$ do not depend on the particular representation of $X^N$.  By 
analogy with the Karhunen--Lo\`eve decomposition of Gaussian processes (see e.g. p57 in~\cite{hida93gaussian}), we call the functions 
$t\mapsto \phi_t(i)$ the {\it discrete principal components} of $X$. 

\begin{Prop}\label{p:principal-comp}
The finite or countable collection of functions $\{t\mapsto \phi_t(i),\ i \in S_N,\ t\in T\},\ N\in \bbN\cup\{\infty\}$
can be the discrete principal components of an $\alpha$--Fr\'echet process, if and only if, the representation
$\{\phi_t\}_{t\in T}\subset L_+^\alpha(S_N,\lambda_N)$ is minimal.
\end{Prop}

\noi The proof is trivial. We state this result to emphasize that not every collection of non--negative
functions can serve as discrete principal components. The {\it minimality} constraint can be viewed
as the counterpart of the {\it orthogonality} condition on the principal components in the Gaussian case. 
The following two examples illustrate typical spectrally discrete and spectrally continuous processes.

\begin{Example} Let $Z_i,\ i\in \bbN$ be independent standard $\alpha$--Fr\'echet variables and let $g_t(i)\ge 0,\ t\in T$ be such that
$\sum_{i\in \bbN} g_t^\alpha(i) <\infty$, for all $t\in T$.  It is easy to see that the $\alpha$--Fr\'echet process
$$
 X_t := \bigvee_{i\in\bbN} g_t(i) Z_i \equiv \Eint{\bbN} g_t dM_\alpha,\ \ t\in T,
$$
is {\it spectrally discrete}.  That is, $X = \{X_t\}_{t\in T}$ has trivial spectrally continuous component.  Indeed, this follows from Theorem \ref{thm:relation} since the mapping $\Phi$ therein is onto, and thus the set $\Phi(\bbN) = S_{I,N}$ is necessarily countable.
\end{Example}

\begin{Example} Consider the well--known $\alpha$--Fr\'echet {\it extremal process} ($\alpha>0$):
\equh\label{rep:extremalProcess}
\{X_t\}_{t\in\mathbb R_+} \defd {\Big\{}\Eintt_{\mathbb R_+}\ind_{(0,t]}(u) M_\alpha(du) {\Big\}}_{t\in \mathbb R_+}\,,
\eque
where $M_\alpha$ has the Lebesgue control measure on $\mathbb R_+$. The process $X=\{X_t\}_{t\in\mathbb R_+}$ can
be viewed as the max--stable counterpart to a sum--stable L\'evy process. This is because $X$ has \textit{independent max--increments}, 
i.e., for any $0=t_0<t_1<\dots<t_n$, 
\[
 (X_{t_1},\dots,X_{t_n}) \eqd (\xi_1,\xi_1\vee\xi_2,\dots,\xi_1\vee\dots\vee\xi_n)\,,
\]
where $\xi_i = M_\alpha((t_{(i-1)},t_i])$,\ $1\le i\le n$.
The representation in~\eqref{rep:extremalProcess} is minimal but its support is not standardized.  Let
$$
 f_t(s) := s^{-1/\alpha} \ind_{(0,t]}(\log(1/s)),\ \ s\in (0,1), 
$$
and observe that $f_t(s) \in L_+^\alpha((0,1),ds)$.  By using a change of variables one can show that
$$
 \{X_t\}_{t\in\mathbb R_+} \defd{\Big\{} \Eintt_{(0,1)} f_t(s) M_\alpha(ds){\Big\}}_{t\in \mathbb R_+},
$$
where the last representation is minimal and has standardized support.  Thus, the
$\alpha$--Fr\'echet extremal process $X$ is {\it spectrally continuous}.
\end{Example}
\subsection{Classification via co--spectral functions}\label{sec:cospectral}

Here we present a characterization of $\alpha$--Fr\'echet processes
based on a different point of view. Namely, instead of focusing on the
spectral functions $s\mapsto f_t(s)$, we now consider the
\textit{co--spectral functions} $t\mapsto f_t(s)$, which are functions of $t$,
with $s$ fixed. To be able to handle the co--spectral functions, we suppose that $T$ is a {\it separable} metric
space with respect to a metric $\rho_T$ and let ${\cal T}$ be its Borel $\sigma$--algebra.  
We say that the spectral representation $\{f_t(s)\}_{t\in T}  \subset L_+^\alpha(S,\mu)$
is jointly measurable if the mapping $(t,s)\mapsto f_t(s)$ is measurable w.r.t.\ the
product $\sigma$--algebra $\calT\otimes{\cal S} := \sigma( {\cal T}\times {\cal S})$.  
The following result clarifies the connection between the joint measurability of the spectral functions
$f_t(s)$ and the measurability of its corresponding $\alpha$--Fr\'echet process.

\begin{Prop}\label{p:measurability}  Let $(S,\mu)$ be a standard Lebesgue space and $M_\alpha$  ($\alpha>0$) 
 be an $\alpha$--Fr\'echet random sup--measure on $S$ with control measure $\mu$.  As above, let $(T,\rho_T)$ be a separable metric space. \\ 
\itemnumber i Let $X=\indt X$ have a spectral representation $\indt f\subset\laps$ as in~\eqref{rep:extremalRep}. Then, $X$ has a measurable modification if and only if $\{f_t(s)\}_{t\in T}$ has a 
jointly measurable modification, i.e., there exists a $\calT\otimes\calB_S-$measurable mapping  $(s,t)\mapsto g_t(s)$, 
such that $f_t(s) = g_t(s)$ $\mu\ae$ for all $t\in T$.  

\itemnumber  {ii} If an $\alpha$--Fr\'echet process $X=\{X_t\}_{t\in T}$ has a measurable modification, then it satisfies Condition S (see Definition~\ref{d:Cond-S}), and
hence it has a representation as in \eqref{rep:extremalRep}.
\end{Prop}

\noindent The proof is given in Appendix.
The above result shows that for a measurable $\alpha$--Fr\'echet process $X=\{X_t\}_{t\in T}$, one 
can always have a representation as in \eqref{rep:extremalRep}, with jointly measurable spectral representations. Conversely, any $X$ as in \eqref{rep:extremalRep}
with measurable spectral functions has a measurable modification.

Let now $\lambda$ be a $\sigma$--finite Borel measure on $T$. We will view each
$f_\cdot(s)$ as an element of the classes $\lopT$ of non--negative $\calT$--measurable functions,
identified with respect to equality $\lambda$--almost everywhere.  Recall that a set
$\calP\subset\lopT$ is said to be a \textit{positive cone} in $\lopT$, if $c\calP \subset \calP$ 
for all $c\ge 0$. Two cones $\calP_1$ and $\calP_2$ are \textit{disjoint} if $\calP_1\cap\calP_2 = \{0\}$.

We propose a general strategy for classification of $\alpha$--Fr\'echet
 processes, based on any collection of disjoint positive cones 
$\calP_j\subset\lopT,\ 1\le j\le n$. 
For any $\alpha$--Fr\'echet process $X = 
 \{X_t\}_{t\in T}$ with jointly measurable representation of full support $\{f_t(s)\}_{t\in T} \subset \laps$, we say the representation has a {\it co--spectral decomposition} w.r.t.\ $\{\calP_j\}_{1\leq j\leq n}$,
if there exist measurable sets $S\topp j, 1\leq j\leq n$, such that
\equh\label{eq:Sis}
 S\topp j \subset \{s\in S\, :\,  f_.(s) \in {\cal P}_j \},\ \ 1\le j \le n\qmand
\mu\Big(S\setminus\bigcup_{j=1}^nS\topp j\Big) = 0\,.
\eque
The sets $S\topp j,1\leq j\leq n$ are modulo $\mu$ disjoint. Indeed,
Let $A:= \{ s\in S\, :\, f_\cdot(s) \equiv 0\}$ and note that $\mu(A) = 0$ by the fact $\supp\{f_t,\ t\in T\} = S$ modulo $\mu$ and Fubini's Theorem.
Since ${\cal P}_j \cap {\cal P}_k = \{0\}$, we have that $S\topp j\cap S\topp k = A$ for
all $1\le j \not = k \le n$.  That is, the space $S$ is partitioned into $n$ modulo $\mu$
disjoint components:
\begin{equation}\label{e:S-co-spec}
 S = S\topp1 \cup \cdots \cup S\topp n\mod\mu,\ \mbox{ with }\ \mu(S\topp j \cap S\topp k) = 0,\ j\not = k.
\end{equation}
This yields the {\it decomposition}:
\begin{equation}\label{e:X-co-spec}
\indt X \eqd \bccbb{X_t^{(1)} \vee \cdots \vee X_t^{(n)}}_{t\in T}\,, 
\end{equation}
with:
\[
X_t^{(j)} := \Eint{S\topp j} f_t (s)M_\alpha(ds), \ \ 1\le j\le n\,,\ \ \forall t\in T.
\]
Note that given a spectral representation $\indt f\subset\lap\smu$, the co--spectral decomposition is defined modulo $\mu$--null sets and the induced decomposition is invariant w.r.t.\ the versions of the decomposition. Namely, if there is another co--spectral decomposition w.r.t.~$\{\calP_j\}_{1\leq j\leq n}$, say $S = \bigcup_{1\leq j\leq n}\widetilde S\topp j\mod\mu$, then from~\eqref{eq:Sis} and the disjointness of $\{\calP_j\}_{1\leq j\leq n}$, it follows that $\mu(\widetilde S\topp j\cap S\topp j) = 0, 1\leq j\leq n$. This yields the same decomposition~\eqref{e:S-co-spec}. 

Moreover, the decomposition is invariant w.r.t.\ the choice of spectral representation.
\begin{Thm}\label{thm:cospectralDecomp}
Suppose $\{\calP_j\}_{1\leq j\leq n}$ are disjoint positive cones in $\lopT$.
For any $\alpha$--Fr\'echet process $\indt X$ with measurable spectral representation $\indt f\subset\lap\smu$, suppose $\indt f$ has a co--spectral decomposition w.r.t.\ $\{\calP_j\}_{1\leq j\leq n}$. Then,

\itemnumber i the decomposition~\eqref{e:X-co-spec} is unique in distribution. 

\itemnumber {ii} the components $\indt {X\topp j}, 1\leq j\leq n$ are independent $\alpha$--Fr\'echet processes.
\end{Thm} 
The proof is given in Appendix. In the special case when $n=1$, Theorem~\ref{thm:cospectralDecomp} yields the following:
\begin{Coro}\label{coro:cospectral}
Let $X=\indt X$ be an $\alpha$--Fr\'echet process with two
jointly measurable representations  $\{f_t^{(i)}(s)\}_{t\in T} \subset\lap(S_i,\mu_i)$, $i=1,2$.
Consider a positive cone ${\cal P}\subset \lopT$.
If $f_\cdot^{(1)}(s) \in {\cal P}$, for $\mu_1$--almost all $s\in S_1$, then
$f_\cdot^{(2)}(s) \in {\cal P}$, for $\mu_2$--almost all $s\in S_2$.
\end{Coro}
Corollary \ref{coro:cospectral} can be used to distinguish between various $\alpha$--Fr\'echet processes in terms of 
their co--spectral functions.  For example, any measurable representation of the $\alpha$--Fr\'echet {\it extremal
process} in \eqref{rep:extremalProcess} should involve simple indicator--type co--spectral functions with one jump down to zero. 
The next result shows another application of Corollary~\ref{coro:cospectral}. 

\begin{Coro} Consider the moving maxima $\alpha$--Fr\'echet random fields:
$$
\{X_t\}_{t\in\bbR^d} \eqd {\Big\{}\eint{\bbR^d} f(t-s) M_\alpha(ds){\Big\}}_{t\in\bbR^d} 
\ \ \mbox{ and } \ \ \{Y_t\}_{t\in \bbR^d}
\eqd {\Big\{} \eint{\bbR^d} g(t-s)M_\alpha(ds){\Big\}}_{t\in\bbR^d},
$$
with $d\in\bbN$, where $f$ and $g$ belong to  $L_+^\alpha(\bbR^d,\lambda)$.  
Here $M_\alpha$ is a an $\alpha$--Fr\'echet random sup--measure on $\bbR^d$ with the 
Lebesgue control measure.
We have $\{X_t\}_{t\in T} \eqd \{Y_t\}_{t\in T}$, if and only if
$g(x) = f(x+\tau)$, almost all $x\in \bbR^d$, with some fixed $\tau\in \bbR^d$.
\end{Coro} 
\begin{proof}
The `if' part is trivial.  To prove the `only if' part, introduce the cone 
${\cal P}_f =\{ cf(\cdot + \tau),\ c\ge 0,\ \tau\in\bbR^d \}$.
Corollary \ref{coro:cospectral} implies that $g(\cdot) \in {\cal P}_f$, and hence $g(x)= c f(x + \tau)$.  Since 
$$
\|X_0\|_\alpha^\alpha = 
 \int_{\bbR^d} g^\alpha (x) dx = \int_{\bbR^d} f^\alpha(x) dx,
$$ 
it follows that $c=1$. This completes the proof.
\end{proof}

Theorem~\ref{thm:cospectralDecomp} is a general result in the sense that the cones $\{\calP_j\}_{1\le j \le n}$ may be associated with various
properties of the co--spectral functions $t\mapsto f_t(s)$ of the process $X$.  If $T\equiv \bbR^d,\ d\ge 1$, for example, one can consider
the cones of co--spectral functions that are: {\it differentiable}, {\it continuous}, {\it integrable}, or {\it $\beta$--H\"older continuous}.
Every choice of cones leads to different types of classifications for measurable $\alpha$--Fr\'echet processes or fields $X=\{X_t\}_{t\in T}$.  
We conclude this section by giving two important examples of classifications, motivated by existing results in the literature on
sum--stable processes.
\begin{Rem}
Note that, instead of~\eqref{e:S-co-spec}, one may want to define $S\topp j\defe\{s:f_\cdot(s)\in\calP_j\},1\leq j\leq n$. However, for certain cones, the $S\topp j$'s defined in this way may not be measurable. See Example~\ref{sec:pos-null}.
\end{Rem}
\begin{Example}[\sc Conservative--dissipative decomposition]\label{sec:cons-diss}
Let $X = \{X_t\}_{t\in T}$ be an
$\alpha$--Fr\'echet process with measurable representation 
$\{f_t(s)\}_{t\in T}\subset L_+^\alpha(S,\mu)$. Consider the following partition of 
the set $S = C\cup D$ with
\begin{equation}
C  \defe  {\Big\{}s: s\in S\,,\int_Tf_t^\alpha(s)\lambda(dt) = \infty {\Big\}}\label{decomp:C,D}\ \ \ 
\mbox{ and } \ \ \ D  \defe  {\Big\{}s: s\in S\,,\int_Tf_t^\alpha(s)\lambda(dt) < \infty {\Big\}}\,.
\end{equation}
Note that $C$ and $D = S\setminus C$ are both $\calS$--measurable since $f_t(s)$ is jointly measurable.  Observe that this partition of $S$ yields the decomposition:
\equh\label{decomp:CD}
 \indt X\eqd \left\{X^C_t\vee X^D_t\right\}_{t\in T}\,,
\eque
where $X^C = \indt{X^C}$ and $X^D = \indt{X^D}$ are defined as:
\equh\label{decomp:CD2}
X_t^C = \Eintt_Cf_tdM_\alpha\ \quad \mbox{ and }\ \quad X_t^D = \Eintt_Df_tdM_\alpha,\ \forall t\in T.
\eque
Here $M_\alpha$ is an $\alpha$--Fr\'echet random sup-measure with control measure $\mu$. 

The decomposition in \eqref{decomp:CD} corresponds to the general decomposition
in \eqref{e:X-co-spec}. Indeed, the co--spectral functions of the component 
$X^D$ belong to the positive cone of {\it integrable} functions, while those of $X^C$ belong to
the cone of non--integrable functions. By Theorem~\ref{thm:cospectralDecomp}, the decomposition \eqref{decomp:CD} does
not depend on the choice of the representation.  The components $X^C$ and $X^D$ of $X$ 
are independent and they are called the {\it conservative} and {\it dissipative} parts of $X$,
respectively. The Decomposition~\eqref{decomp:CD} is referred to as the {\it conservative--dissipative} decomposition.
\end{Example}
\begin{Example}[\sc Positive--null decomposition]\label{sec:pos-null}
Following Samorodnitsky \cite{samorodnitsky05null}, consider $T = \mathbb R$ or $\mathbb Z$. Introduce the class ${\cal W}$
of {\it positive} weight functions $w:T\to\bbR_+$:
\begin{equation}\label{e:W-def}
 {\cal W} := {\Big\{}w: \int_{T} w(t) \lambda(dt) = \infty,\ w(t)\mbox{ and }
 w(-t)\mbox{ are non--decreasing on $T\cap (0,\infty)$}{\Big\}}.
\end{equation}
Now we consider the cone
$$
 {\cal P}_{\rm pos} :=\Big\{ f \in L_+^0(T,\lambda)\, :\, \int_T w(t) f_t^\alpha \lambda(dt) =\infty,
 \mbox{ for all } w\in {\cal W} \Big\}
$$
and its complement cone
${\cal P}_{\rm null} := \{0\}\cup (L_+^0 (T,\lambda) \setminus {\cal P}_{\rm pos})$.

This choice of cones yields the decomposition
\begin{equation}\label{e:pos-null}
\{X_t\}_{t\in T} \eqd\{X_t^{\rm pos} \vee X_t^{\rm null}\}_{t\in T},
\end{equation}
where
\begin{equation}\label{e:pos-null-1}
 X_t^{\rm pos} := \Eint{P} f_t(s) M_\alpha(ds)\ \ \mbox{ and } \ \
 X_t^{\rm null} := \Eint{N} f_t(s) M_\alpha(ds),\ \forall t\in T\,,
\end{equation}
with $P$ and $N$, measurable subsets of $S$, satisfying $\mu(P\cap N) = 0$, $\mu(S\setminus (P\cup N)) = 0$ and 
\begin{equation}\label{e:P-N}
f_\cdot(s)\in\calP_{\rm pos},\forall s\in P\qmand f_\cdot(s)\in \calP_{\rm null}, \forall s\in N\,.
\end{equation}
The components $X^{\rm pos}=\{X_t^{\rm pos}\}_{t\in T}$ and
$X^{\rm null}=\{X_t^{\rm null}\}_{t\in T}$ in \eqref{e:pos-null-1} are said to be the {\it
positive} and {\it null} components of the process $X$, respectively. By Theorem~\ref{thm:cospectralDecomp}, Decomposition~\eqref{e:pos-null} does not depend on the choice of the measurable representation $\{f_t(s)\}_{t\in T}
\subset L_+^\alpha(S,\mu)$. It is 
referred to as the {\it positive--null} decomposition. 

Note that, a technical difference between this example and Example~\ref{sec:cons-diss} is that the set $\widetilde P\defe\{s:f_\cdot(s)\in\calP_{\rm pos}\}$ may not be measurable, even when $f_t(s)$ is jointly measurable.
\end{Example}

In the following section, we will study the above decompositions in more detail, for the case of stationary max--stable processes.
\section{Classification of Stationary $\alpha$--Fr\'echet Processes}\label{sec:stationary}
In this section, we focus on stationary, measurable max--stable processes $X = \{X_t\}_{t\in T}$, 
where $T=\bbR$ or $T=\bbZ$ is equipped with the Lebesgue or the counting measure $\lambda$, 
respectively. In this case, the process $X$ can be associated with a non--singular flow. 
Therefore, as in the symmetric $\alpha$--stable case, the ergodic theoretic properties of the flow
yield illuminating structural results.

\subsection{Non--singular flows associated with max--stable processes}
\label{sec:flows}

Following Rosi\'nski \cite{rosinski95structure} (see also Appendix A in~\cite{pipiras04stable}), we recall some notions from ergodic theory.
\begin{Def}\label{def:flow}
A family of functions $\phi=\{\phi_t\}_{t\in T}$, $\phi_t:S\to S$ for all $t\in T$, is a flow on $(S,\calB,\mu)$ if\\
\itemnumber i $\phi_{t_1+t_2}(s) = \phi_{t_2}(\phi_{t_1}(s))\,, \forall t_1,t_2\in T\,,s\in S$.\\
\itemnumber {ii} $\phi_0(s) = s\,, \forall s\in S$.\\
A flow $\phi$ is said to be \textit{measurable} if $\phi_t(s)$ is a measurable map from $T\times S$ to $S$;
A flow $\phi$ is said to be \textit{non--singular} if $\mu(\phi_t\inv(A)) = 0 \Leftrightarrow 
\mu(A) = 0, \forall A\in \calB\,,t\in T$.
\end{Def}

\noindent
The next result relates the spectral functions of stationary $\alpha$--Fr\'echet processes to 
flows.
\begin{Thm}\label{thm:flow}
Let $\pxt$ be a stationary $\alpha$--Fr\'echet process.  Suppose that $X$ has a measurable 
representation $\fft\subset \lap(\ssS,\ssL)$, which is minimal, with standardized support.
Then, there exist a unique, modulo $\ssL$, non--singular and measurable 
flow $\{\phi_t\}_{t\in T}$ such that for each $t\in T$,
\equh\label{rep:flowRepSS1}
 f_t(s) = \left(\dfrac\ssL{\phi_t}\ssL\right)^{1/\alpha}(s)(f_0\circ\phi_t)(s)\,,\quad \ssL\ae\,.
\eque
\end{Thm}
Theorem~\ref{thm:flow} is stronger than Theorem~6.1 in~\cite{dehaan86stationary}, where the measurability is not considered and the flow structure is not explicitly explored.
The proof is given in Appendix~\ref{sec:proofStationary}.  For the readers familiar with Rosi\'nski's work~\cite{rosinski95structure}, this result is similar to Theorem 3.1 therein. 
In view of this result, we will say that a stationary $\alpha$--Fr\'echet measurable
process $\indt X$ is {\it generated} by the non--singular measurable flow $\indt\phi$ on
$(S,\mu)$ if it has a spectral representation $\{f_t\}_{t\in T}\subset L_+^\alpha(S,\mu)$, where:
\equh\label{rep:flowRep1}
 f_t = \left(\dfrac\mu{\phi_t}\mu\right)\tia(f_0\circ\phi_t),\quad\mu\ae,
\eque
and 
\equh\label{rep:flowRep2}
\supp\{f_0\circ\phi_t:t\in T\} = S,\quad\mu\ae
\eque
Note that in the representation~\eqref{rep:flowRep1} and~\eqref{rep:flowRep2}, we do not assume 
$\indt f$ to be minimal. However, the minimality plays a crucial role in the proof of the
existence of flow representations in Theorem~\ref{thm:flow}.

\begin{Def}\label{def:equivalenceOfFlows} We say two measurable non--singular 
flows $\indt{\phi\topp1}$ and $\indt{\phi\topp2}$ on $(S_i,\mu_i), i = 1,2$, are equivalent, 
written $\indt{\phi\topp1}\sim^\Phi \indt{\phi\topp2}$, if there exists a measurable map 
$\Phi:S_2\to S_1$ such that:

\itemnumber {i} There exist $N_i\subset S_i$ with $\mu_i(N_i) = 0, i=1,2$ such that $\Phi$ is a
Borel isomorphism between $S_2\setminus N_2$ and $S_1\setminus N_1$.\\
\itemnumber {ii} $\mu_1$ and $\mu_2\circ\Phi\inv$ are mutually absolutely continuous.\\
\itemnumber {iii} $\phi_t\topp1\circ\Phi = \Phi\circ\phi_t\topp2 \mu_2\ae$ for each $t\in T$.
\end{Def}

The next result shows the connection between different flows generating the same stationary 
$\alpha$--Fr\'echet process $\indt X$. The proof is given in Appendix~\ref{sec:proofStationary}.
 
\begin{Prop}\label{prop:equiFlow}
Let $\indt X$ be a measurable stationary $\alpha$--Fr\'echet process.\\
\itemnumber i Suppose $\indt{\phi\topp 1}$ is a flow on $(S_1,\mu_1)$ and $\indt X$ is generated by 
$\indt{\phi\topp 1}$ with spectral function $f_0\topp1\in\lap(S_1,\mu_1)$. If $\indt{\phi\topp2}$ is 
another flow on $(S_2,\mu_2)$ and it is equivalent to $\indt{\phi\topp1}$ via $\Phi$, then $\indt X$ 
can also be generated by $\indt{\phi\topp2}$ with the spectral function
\equh\label{eq:equivalence2}
f_0\topp2(s) = \left(\dfrac{\mu_1}{\Phi}{\mu_2}(s)\right)^{1/\alpha}\left(f_0\topp1\circ\Phi\right)(s)\,.
\eque
Moreover, if $\indt {f\topp1}$ is minimal, then $\indt {f\topp2}$ is minimal.\\
\itemnumber {ii} If $\indt X$ has two measurable minimal representations generated by flows $\indt{\phi\topp i}$ on $(S_i,\mu_i)$ for $i=1,2$, then $\indt{\phi\topp1}\sim^\Phi\indt{\phi\topp2}$ and~\eqref{eq:equivalence2} holds, for some $\Phi$ satisfying conditions in Definition~\ref{def:equivalenceOfFlows}.
\end{Prop}

\begin{Rem} Not all flow representations are minimal.  Proposition \ref{prop:equiFlow} 
shows, however, that any two flows corresponding to minimal representations of the same 
$\alpha$--Fr\'echet process are equivalent. \end{Rem}


\subsection{Decompositions induced by non--singular flows}\label{sec:decompositions}

The decompositions introduced in Examples \ref{sec:cons-diss} and \ref{sec:pos-null}
are motivated by corresponding notions from ergodic theory. 

\begin{Def}\label{def:wandering} Consider a measure space $(S,\mu)$ and a 
measurable, non--singular map $\phi : S \to S$.  A measurable set 
$B\subset S$ is said to be:

\itemnumber i {\it wandering}: if $\phi^{-n} (B),\ n=0,1,2,\cdots$ are disjoint.

\itemnumber {ii} {\it weakly wandering}: if $\phi^{-n_k} (B),\ n_k\in \bbN$ 
are disjoint, for an infinite sequence $0=n_0 < n_1 <\cdots$.
\end{Def}
Now we give two decompositions for max--stable processes. Their counterparts for sum--stable processes have been thoroughly studied (see~\cite{rosinski95structure} and~\cite{samorodnitsky05null}).\medskip

{\sc Hopf (conservative--dissipative) decomposition.}
The map $\phi$ is said to be {\it conservative} 
if there is no {\it wandering} measurable set $B\subset S$, with positive measure $\mu(B)>0$.
One can show that for any measurable, non--singular map $\phi:S\to S$, there exists a 
partition of $S$ into two disjoint measurable sets $S = C\cup D$, $C\cap D = \emptyset$ such that:
{\it (i)} $C$ and $D$ are $\phi-$invariant; {\it (ii)} $\phi:C\to C$ is conservative and 
$D = \cup_{k\in \bbZ} \phi^k(B),$ for some wandering set $B\subset S$.
This decomposition is unique (mod $\mu$) and is called the {\it Hopf decomposition} of 
$S$ with respect to $\phi$.  If the component $C$ is trivial, i.e.\ $\mu(C)=0$, then $\phi$
is said to be {\it dissipative}.  The restrictions $\phi:C\to C$ and $\phi:D\to D$ are the
{\it conservative} and {\it dissipative} components of the mapping $\phi$, respectively.

Now, given a jointly measurable, non--singular flow $(t,s)\mapsto \phi_t(s),\ t\in T,\ s\in S$,
one can consider the Hopf decompositions $S=C_t\cup D_t$ for each $\phi_t,\ t\in T\setminus\{0\}$.
By the measurability however, it follows that $\mu(C_t \Delta C) = \mu(D_t \Delta D) = 0$, for some
$C\cap D = \emptyset$, $S = C\cup D$ (see e.g.\ \cite{rosinski95structure,Krengel85ergodic}).
One thus obtains that any measurable non--singular flow $\{\phi_t\}_{t\in T}$ has a Hopf 
decomposition $S= C\cup D$, where $\phi^C:= \{\phi_t\vert_C\}_{t\in T}$ and $\phi^D:= 
\{\phi_t\vert_D\}_{t\in T}$ are {\it conservative} and {\it dissipative} flows, respectively.

The following result is an immediate consequence from the 
proofs of Theorem 4.1 and Corollary 4.2 in Rosi\'nski \cite{rosinski95structure}.

\begin{Thm}\label{thm:cons-diss}
Let $X=\{X_t\}_{t\in T}$ be a stationary $\alpha$--Fr\'echet process
with measurable representation $\{f_t(s)\}_{t\in T}\subset L_+^\alpha(S,\mu)$ of full support.
Then: 

\itemnumber {i} $X$ is generated by a conservative flow, if and only if,
$$
 \int_T f^\alpha_t(s) \lambda(dt) = \infty,\ \ \mbox{ for $\mu$--almost all $s\in S$; }
$$
\itemnumber {ii} $X$ is generated by a dissipative flow, if and only if,
$$
 \int_T f^\alpha_t(s) \lambda(dt) < \infty,\ \ \mbox{ for $\mu$--almost all $s\in S$. }
$$
\itemnumber {iii} If $X$ is generated by a conservative (dissipative) flow in one representation, 
then so is the case for any other measurable representation of $X$.
\end{Thm}

This result justifies the terminology in the {\it conservative--dissipative} decomposition of
Example~\ref{sec:cons-diss}.  In particular, the sets $C$ and $D$ in \eqref{decomp:CD} 
correspond precisely to the conservative and dissipative parts in the Hopf decomposition of the
flow $\{\phi_t\}_{t\in T}$ associated with the process $X$. \medskip

{\sc Positive--null decomposition.} Recall the notion of {\it weakly wandering}
set (Definition \ref{def:wandering}).  If one replaces `wandering' by `weakly wandering' in 
the Hopf decomposition, one obtains the so--called {\it positive--null decomposition} of $S$. 
Alternatively, the map $\phi$ is said to be  {\it positive}, if there exists a finite measure 
$\nu \sim \mu$, such that $\phi$ is $\nu$--invariant. In this case, there are no weakly 
wandering sets $B$ of positive $\mu$--measure (or equivalently, $\nu$--measure).  For any 
non--singular map $\phi$, there exists a partition $S = P\cup N$, unique modulo $\mu$, such 
that $P$ and $N$ are disjoint, measurable and $\phi$--invariant.  Furthermore, $\phi:P\to P$ is 
positive, and $N = \cup_{k\ge 0} \phi^{-n_k}(B)$,
for some disjoint $\phi^{-n_k}(B)$'s, where $B$ is weakly wandering. 
The set $N$ ($P$ resp.) is called the null--recurrent (positive--recurrent) part of $S$, w.r.t. the map $\phi$ (see e.g.\ Section 1.4 in \cite{aaronson97introduction}). 

As in the case of the Hopf decomposition, a jointly measurable, non--singular
flow $\{\phi_t\}_{t\in T}$ gives rise to a {\it positive--null decomposition}: $S = P\cup N,$
where $\mu(P_t\Delta P) = \mu(N_t\Delta N) = 0$, for all $t\in T\setminus\{0\}$, and where
$S = P_t \cup N_t$ is the positive--null decomposition of the map $\phi_t,\ t\in T\setminus\{0\}$
(see e.g.\ \cite{samorodnitsky05null,Krengel85ergodic}).

Theorem 2.1 of Samorodnitsky \cite{samorodnitsky05null} about symmetric $\alpha$--stable processes
applies {\it mutatis mutandis} to the max--stable case:

\begin{Thm}\label{thm:pos-null}
Let $X=\{X_t\}_{t\in T}$ be a stationary $\alpha$--Fr\'echet process
with measurable representation $\{f_t(s)\}_{t\in T}\subset L_+^\alpha(S,\mu)$ of full support. Then: \\
\itemnumber {i} $X$ is generated by a positive flow, if and only if, for all $w\in {\cal W}$,
$$
 \int_T w(t) f^\alpha_t(s) \lambda(dt) = \infty,\ \ \mbox{ for $\mu$--almost all $s\in S$, }
$$
where ${\cal W}$ is as in \eqref{e:W-def}.\\
\itemnumber {ii} $X$ is generated by a null flow, if and only if, for some $w\in {\cal W}$,
$$
 \int_T w(t) f^\alpha_t(s) \lambda(dt) < \infty,\ \ \mbox{ for $\mu$--almost all $s\in S$. }
$$
\itemnumber {iii} If $X$ is generated by a positive (null) flow in one representation, 
then so is the case for any other measurable representation of $X$.
\end{Thm}

As in the Hopf decomposition, Theorem \ref{thm:pos-null} shows that the components
$X^{\rm pos}$ and $X^{\rm null}$ in the decomposition \eqref{e:pos-null} are generated by
positive-- and null--recurrent flows, respectively. This is because the sets $P$ and $N$ in
\eqref{e:P-N} yield the positive--null decomposition of a flow $\{\phi_t\}_{t\in T}$
associated with $X$.

\subsection{Structural results, examples and open questions}

Here, we collect some structural results and observations on the interplay between
the three types of classifications of max--stable processes discussed above.
Namely, {\it (i)} continuous--discrete {\it (ii)} conservative--dissipative and 
{\it (iii)} positive--null.  
 
Theorems~\ref{thm:cons-diss} and~\ref{thm:pos-null} imply that the {\it positive} component of a
max--stable process is {\it conservative} and the {\it dissipative} one is {\it null--recurrent}. 
Thus, for a measurable stationary $\alpha$--Fr\'echet process $\indt X$, we have the
decomposition:
\equh\label{decomp:flow}
\indt X \eqd \bccbb{X^{\rm pos}_t \vee X^{C,\rm null}_t \vee X^D_t}_{t\in T}\,,
\eque
where $X^{C}_t = X^{\rm pos}_t \vee X^{C,\rm null}_t$ and 
$X^{\rm null}_t = X^{C,\rm null}_t\vee X^{D}_t$, $t\in T$.
Here $X^{\rm pos},\ X^{C,\rm null}$ and $X^D$ are independent $\alpha$--Fr\'echet processes.
$X^{\rm pos}$ is positive--recurrent and conservative, $X^D$ is dissipative and null--recurrent, 
and $X^{C,\rm null}$ is conservative and null--recurrent. We will see that the $X^D$ is precisely 
the mixed moving maxima. Moreover, we show that the spectrally discrete component has no
conservative--null component $X^{C,\rm null}$. 
 
The following theorem shows that the {\it purely dissipative} stationary 
$\alpha$--Fr\'echet processes are precisely the \textit{mixed moving maxima}.  

\begin{Thm}\label{thm:mmm}
Let $\pxt$ be a measurable stationary $\alpha$--Fr\'echet process.
This process is generated by a dissipative flow if and only if 
there exist a Borel space $W$, a $\sigma$-finite measure $\nu$ on $W$ and a function 
$g\in\lap(W\times T,\nu\otimes\lambda)$ such that
\[
\pxt \eqd \left\{\Eintt_{W\times T}g(x,t+u)M_\alpha(dx,du)\right\}_{t\in T}\,.
\]
Here $M_\alpha$ is an $\alpha$--Fr\'echet random sup-measure on $W\times T$ with the control measure $\nu\otimes\lambda$ and $\lambda$ is the Lebesgue measure if $T=\mathbb R$ and the counting measure if $T=\mathbb Z$. Moreover, one can always choose $(W,\nu)$ and $g$ such that the representation $g_t(x,u) \defe g(x,t+u)$ is minimal.
\end{Thm}
\begin{proof} Since $g\in\lap(W\times T,\nu\otimes\lambda)$, the Fubini's theorem implies
$\int_{T} g(x,t+u)^\alpha \lambda(dt)<\infty$, for almost all $(x,u)\in W\times T$.  This, in view of \eqref{decomp:C,D} implies that $X$ is dissipative..

The `only if' part follows as in the proof of Theorem 4.4 in Rosi\'nski \cite{rosinski95structure} from the results of Krengel \cite{krengel69}.  
\comment{Namely, for every dissipative flow $\indt\phi$ on $(S,
\mu)$, there exists a finite standard Lebesgue space $(W,\nu)$ such that the 
flow $\indt\phi$ is null isomorphic to a flow $\indt\beta$ defined on 
$(W\times T,\nu\otimes \lambda)$ by
\[
\beta_t (x,u) = (x,t+u)\,,\quad(x,u)\in W\times T\,,t\in T\,.
\]
That is, there exists a non--singular invertible map $\Phi:W\times T\to S$ such that $\Phi\circ\beta_t 
= \phi_t\circ\Phi\,,\nu\ae$ for all $t\in T$. Since $d[(\nu\otimes\lambda)\circ\beta_t]/d(\nu\otimes\lambda) = 
1$ for every $t\in T$, $\indt\phi$ and $\indt\beta$ are equivalent and by 
Proposition~\ref{prop:equiFlow} we obtain the desired result. In effect, let $\indt f$ be the 
spectral function in flow representation~\eqref{rep:flowRep1}. Define
\[
g_0 = \left(\dfrac{\mu}{\Phi}{(\nu\otimes\lambda)}\right)^{1/\alpha}(f_0\circ\Phi)\,,\quad g_t = 
g_0\circ\beta_t\,.
\]
Then, by Proposition \ref{prop:equiFlow}, it follows that $\indt X$ has the representation
\[
\indt X \eqd \left\{\Eintt_{W\times T}g_0\circ\beta_t(x,u) M_\alpha(dx,du)\right\}_{t\in T}\,,
\]
and if $\indt{f}$ is minimal, then so is $\indt g$.}
\end{proof}

\begin{Rem} Theorem \ref{thm:mmm} parallels the fact that the class of stationary and dissipative 
symmetric $\alpha$--stable processes is precisely the class of mixed moving averages
(see Theorem 4.4 in \cite{rosinski95structure}).  Recently, Kabluchko \cite{kabluchko08spectral}
established the same result as in Theorem \ref{thm:mmm} by using an interesting {\it association device} between $\alpha$--Fr\'echet ($\alpha\in
(0,2)$) and symmetric $\alpha$--stable processes. 
\end{Rem}

As shown in \cite{stoev08ergodicity}, the mixed moving maxima processes are mixing and
 hence ergodic.  Thus, Theorem \ref{thm:mmm} implies that the dissipative component of a
  max--stable process is mixing. On the other hand, Samorodnitsky \cite{samorodnitsky05null} has shown (Theorem 3.1 therein) that stationary symmetric 
$\alpha$--stable processes are ergodic {\it if and only if} they are generated by a 
null--recurrent flow.  Kabluchko \cite{kabluchko08spectral} (Theorem 8 therein)
has shown that this continues to be the case for stationary $\alpha$--Fr\'echet processes. 

The previous discussion shows that the ergodic and mixing properties of the null and dissipative components are in line with the
decomposition $X_t^{\rm null} = X_t^{D}\vee X_t^{C,\rm null},\ t\in T$.  An example of 
conservative--null flow can be found in \cite{samorodnitsky05null}.  This yields non--trivial
examples of sum-- and max--stable processes that are conservative and null.  We are not aware, 
however, of an example of an ergodic max--stable process that is not mixing. 

The next two results clarify the structure of the stationary {\it spectrally discrete} 
processes in discrete $(T=\bbZ)$ and continuous ($T=\bbR$) time, respectively. 
We first show that for {\it spectrally discrete} stationary max--stable time series, the
conservative--dissipative and positive--null decompositions coincide. That is, such processes
have no {\it conservative--null} components. Moreover, the {\it dissipative} (equivalently 
{\it null--recurrent}) component does not exist if the time series has only {\it finite} 
number of principal components. 

\begin{Prop} \label{p:disc-cons} Let $X = \{X_t\}_{t\in T}$, with $T=\bbZ$ be a
 stationary $\alpha$--Fr\'echet process (time series). \medskip
 
 \itemnumber {i} $X^N$ has no conservative--null component, i.e.\ $X^{N,C,\rm null} = 0$.\\
 \itemnumber {ii} If $1\leq N<\infty$, then $X^N$ is necessarily conservative, and equivalently, positive recurrent.
\end{Prop}
\begin{proof} Without loss of generality, suppose $X = X^N$ and 
let $\{f_t(s)\}_{t\in T} \subset L_+^\alpha(S_N,\lambda_N)$ be a minimal
representation with standardized support for $X$.  We have that
$$
 f_t = \left(\dfrac{\lambda_N}{\phi_t}{\lambda_N} \right)^{1/\alpha} f_0\circ\phi_t,
$$
where $\phi_t:S_N\to S_N$ is a non--singular flow on $(S_N,\lambda_N)$.  Since 
$S_N\subset\mathbb N$ and $\lambda_N$ is the counting measure,
the non--singular transformations are necessarily measure--preserving, i.e., permutations. 
Thus the term $d(\lambda_N\circ \phi_t)/d\lambda_N \equiv  1$ and $f_t(s) = f_0\circ\phi_t(s)$. 

We start by proving {\it (ii)}. Since $\phi_1:\{1,\cdots,N\} \to \{1,\cdots,N\}$ is a permutation, it has a finite invariant measure and hence the flow $\indt\phi$ is positive--recurrent and hence conservative. 

Now we prove {\it (i)}. Note that when $1\leq N<\infty$, we have shown in {\it (ii)} that $X^N$ is conservative and positive--recurrent. For $N = \infty$, we consider two cases. First we suppose that for every $s\in S_N$, the recurrent time 
\equh\label{eq:taus}
\tau_s\defe\inf \{t>0:\phi_t(s) = s\}
\eque
is finite. Let $\mathfrak O(s)$ denote the orbit of state $s$ w.r.t. flow $\indt\phi$, i.e., $\mathfrak O(s)\defe\{\phi_t(s):t\in T\}$. 
Every orbit of $\indt\phi$ is $\tau_s$--periodic, i.e., $|\mathfrak O(s)| <\infty$. Since $N = \infty$, the total number of different orbits must be infinite. Enumerate all the orbits by $\mathfrak O_1,\mathfrak O_2,\cdots$, so that $\mathfrak O(s) = \mathfrak O_{\pi(s)}$ with $\pi: S_N\to\mathbb N$ and $S_N = \bigcup_{k\in\mathbb N}\mathfrak O_k$. Observe that the orbits are disjoint. We now define a finite invariant 
measure on $S_N$, equivalent to the counting measure: 
\[
\widetilde\lambda(\{s\})\defe 2^{-\pi(s)}\frac1{|\mathfrak O_{\pi(s)}|}\,,\forall s\in S_N\,.
\]
This measure is clearly is invariant on each $\mathfrak O_k$, for all $k\in\mathbb N$.  Since
$\widetilde\lambda(\mathfrak O_k) = 2^{-k}$, the measure $\widetilde \lambda$ is finite and it is
clearly equivalent to the counting measure.  Thus, $X^N$ is positive and conservative.

On the other hand, suppose that there exists a state $s$ with $\tau_s=\infty$. Then,
its orbit is infinite and non--recurrent., i.e., $|\mathfrak O_k(s)| = \infty$. Then, 
the flow $\indt\phi$ is both null--recurrent and dissipative on $\mathfrak O_k(s)$. Indeed, the null recurrence
follows from the fact that there is no positive finite invariant measure on $\mathfrak O_k(s)$. 
The dissipativity follows from the remark that $\mathfrak O_k(s) = \bigcup_{j\in\mathbb
 Z}\phi_j(s)$ is a disjoint union. We have thus shown that $\indt\phi$ is dissipative and null--recurrent on non--recurrent orbits. 
\end{proof}

The following result shows that the {\it continuous--time} stationary, measurable and {\it spectrally discrete} max--stable processes 
are trivial.

\begin{Thm}\label{thm:continuousTime}
Let $X = \{X_t\}_{t\in T}$, with $T = \bbR$ be a
stationary and measurable $\alpha$--Fr\'echet process.  If $N\geq 1$, then it must be $N = 1$. That is, the spectrally discrete component $X^N$
is the random constant process:
  $\{X^N_t\}_{t\in \bbR} \eqd \{Z\}_{t\in\bbR}$, for some $\alpha$--Fr\'echet
  variable $Z$.
\end{Thm}
\begin{proof}
Let $\indt f$ and $\indt\phi$ be as in Proposition~\ref{p:disc-cons}. Observe moreover that, in this case, 
the $\phi_t$'s are measure--preserving bijections, and in view of Theorem \ref{thm:flow}, the flow $\{\phi_t(s)\}$ is measurable. For any {\it fixed} $s\in S_N$, consider $\tau_s$ defined in~\eqref{eq:taus}. The proof consists of three steps.

\itemnumber i {\it We show first that $\tau_s = 0$ implies $\phi_t(s) \equiv s$, for all $t\in\mathbb R$.}  Indeed, suppose that $\tau_s = 0$ and note 
that, by definition, for all $n>0$, there exist $0<t_{n,1}<t_{n,2}<1/n$ such that $\phi_{t_{n,1}}(s) = \phi_{t_{n,2}}(s) = s$. Set
$T_0\defe\bigcup_{n\in\mathbb N}\bigcup_{k\in\mathbb Z}\{t_{n,1}+k(t_{n,2}-t_{n,1})\}$. It follows that $T_0$ is dense in $\mathbb R$ and $\phi_t(s) = s$,
for all $t\in T_0$. Hence $f_t(s) = f_0\circ\phi_t(s) = f_0(s)$, for all $t\in T_0$. Now, we define a new $\alpha$--Fr\'echet process $Y = \indt Y$:
\[
\indt Y \eqd\bccbb{\Eintt_{S_N}\ind_{\{\cdot = s\}}\circ\phi_t(r) M_\alpha(dr)}_{t\in T}\,.
\] 
Since $\indt \phi$ is a flow, $\phi_t$ is invertible, for any $t\in T_0$. Hence, for all $t\in T_0$, we have $\phi_t(r) = \phi_t(s)\equiv s$ 
if and only $r = s$.  This shows that, for all $t\in T_0$,
$$
 \ind_{\{\cdot = s\}}\circ\phi_t(r) \equiv \ind_{\{\phi_t(r) = s\}} = \ind_{\{r=s\}} \equiv \ind_{\{\cdot = s\}} \circ \phi_0(r),
$$
which implies that $Y_t = Y_0,$ almost surely, $\forall t\in T_0$. Moreover, as $\indt \phi$ is measurable, so is $\indt Y$ by
Proposition~\ref{p:measurability}. Also, $Y=\indt Y$ is stationary, since it is generated by a measure preserving flow.  Thus, the stationarity and
measurability of $Y$ imply that it is {\it continuous in probability} (see Theorem~3.1 in~\cite{stoev08ergodicity}).  
This, and the fact that $Y_t = Y_0$, a.s., for all $t$ in a dense sub--set $T_0$ of $\mathbb R$, imply  that $Y_t = Y_0$, a.s., for all $t\in\bbR$.
Therefore, for the spectral functions, we obtain $\ind_{\{\phi_t(r) = s\}} = \ind_{\{r = s\}}\,,\forall r\in S_N, t\in\mathbb R$.
This shows that $\phi_t(s) = s,\forall t\in\mathbb R$.

\itemnumber {ii} {\it We show next that $\tau_s>0$ implies $\phi_{\tau_s}(s) = s$.}  
Suppose that $\phi_{\tau_s}(s) \not = s$.  Then, as above, there exist $t_1,t_2\in(\tau_s,\tau_s+\tau_s/2)$ such that $\phi_{t_1}(s) = \phi_{t_2}(s) = s$.
But it follows that $\phi_{t_1+k(t_2-t_1)}(s) = s$ for all $k\in\mathbb Z$.  This, since $\{t_1+k(t_2-t_1)\}_{k\in\bbZ} \cap (0,\tau_s) \not= \emptyset$,
contradicts the definition of $\tau_s$.

\itemnumber {iii} {\it Now, we show that it is impossible to have $\tau_s>0$ for all $s\in S_N$.}
Write $\mathfrak T_s = \{t:\phi_t(s_0) = s,\ \mbox{ for some }s_0\in S_N\}$. Observe that the set $\mathfrak T_s$ is countably infinite for all
$s\in S_N$ such that $\tau_s>0$, since by {\it (ii)} above, $\mathfrak T_s = \{ k\tau_s\}_{k\in\bbZ}$.  
Note also that $\bigcup_{s\in S_N}\mathfrak T_s = \mathbb R$.  However, the assumption that $\tau_s>0$ for all $s\in S_N$ would imply
$\bigcup_{s\in S_N}\mathfrak T_s$ has cardinality of $\mathbb N$ equals that of $\mathbb R$, which is a contradiction.

We now conclude the proof. By {\it (iii)} above, there must exist $s\in S_N$ such that $\tau_{s} = 0$. Set
$\mathfrak R = \{s\in S_N:\phi_t(s) = s,\forall t\in\mathbb R\}$. We have already seen in {\it (i)} that $\tau_{s} = 0$ implies $\phi_t(s) \equiv s$,
for all $t\in\mathbb R$, whence $\mathfrak R$ is $\phi$--invariant. Consider now a new $\alpha$--Fr\'echet process 
\[
 \indt{Y}\eqd\bccbb{\Eintt_{S_N\setminus\mathfrak R}f_t(r)M_\alpha(dr)}_{t\in T}\,.
\]
Since the $f_t$'s, restricted to the $\phi-$invariant set $S_N\setminus\mathfrak R$ yield a minimal representation for $Y= \indt{Y}$ with
standardized support.  This process is generated by the same flow $\{\phi_t\}_{t\in\bbR}$, restricted to $S_N\setminus\mathfrak R$.
Since $\tau_s>0,\ \forall s\in S_N\setminus\mathfrak R$, by {\it (ii)}, it follows that $S_N\setminus\mathfrak R = \emptyset$. 

On the other hand, since $\phi_t(s)\equiv s\,,\forall t\in\mathbb R, s\in\mathfrak R$, the minimality of $\indt f$ implies that
$| \mathfrak R | = |S_N| =1$. Therefore, $\indt X \eqd \{Z\}_{t\in T}$ for some $\alpha$--Fr\'echet random variable $Z$.
\end{proof}

\begin{Example} In contrast with Proposition \ref{p:disc-cons} {\it (i)},
the spectrally discrete component of a stationary $\alpha$--Fr\'echet time series may be dissipative if it 
involves infinite number of principal components. Indeed, by Theorem \ref{thm:mmm}, the moving maxima $X_t:= \Eint{\mathbb Z} f(t+s) M_\alpha(ds) \equiv \bigvee_{i\in\mathbb Z} f(t+i) M_\alpha(\{i\}),$ is dissipative and spectrally discrete, where $M_\alpha$ has the 
counting control measure on $\mathbb Z$.
\end{Example}

\begin{Example} Suppose that $(E,\calE,\mu)$ is a probability space, i.e.\ $\mu(E)=1$. Let
$M_\alpha$ be an $\alpha$--Fr\'echet random sup--measure on $E$ with control 
measure $\mu$, which is defined on a {\it different} probability space. 
Suppose that $\{Y_t\}_{t\in T}$ is a positive stochastic process on $(E,{\cal E},\mu)$ such that
$\E_\mu Y_t^\alpha <\infty,$ for all $t\in T$.  Then, the  $\alpha$--Fr\'echet process:
\equh\label{rep:doublyStochastic}
\indt X \eqd {\Big\{}\Eintt_EY_t(u) M_\alpha(du){\Big\}}_{t\in T}\,,
\eque
is said to be {\it doubly stochastic}.
\end{Example}
One can show that, in~\eqref{rep:doublyStochastic}, if $\indt {Y}$ is stationary, then so is $\indt X$.  The Brown--Resnick processes discussed in next section shows that the converse is not always true. 

\section{Brown--Resnick Processes}\label{sec:BRp}
Consider the following \textit{doubly stochastic process} (see e.g. \cite{kabluchko08stationary} and  \cite{stoev08ergodicity}):
\equh\label{rep:brownResnick}
\indtr X \eqd \left\{\Eintt_Ee^{W_t-\sigma^2_t/2}dM_1\right\}_{t\in \mathbb R}\,.
\eque
Here $W_t$ is a zero--mean Gaussian process defined on the probability space $(E,{\cal E},\mu)$
with variance $\sigma^2_t$. Since $\esp_\mu {\Big(}e^{W_t-\sigma_t^2/2} {\Big)} = 1<\infty$, 
the 1--Fr\'echet process in~\eqref{rep:brownResnick} is well--defined. The processes having 
representation~\eqref{rep:brownResnick} were first introduced by Brown and
Resnick~\cite{brown77extreme} with $W_t$ being the standard Brownian motion. In general, we will 
call $\indtr X$ as in~\eqref{rep:brownResnick} a  \textit{Brown--Resnick 1--Fr\'echet process}.

Kabluchko {\it et al.} \cite{kabluchko08stationary} have shown that if $\indtr W$ has 
\textit{stationary increments}, then the Brown--Resnick process $\indtr X$
in~\eqref{rep:brownResnick} is stationary.  The following interesting result about an arbitrary 
zero--mean Gaussian process with stationary increments and continuous paths is obtained by 
combining the results of \cite{kabluchko08stationary} and our Theorems \ref{thm:cons-diss} and 
\ref{thm:mmm} above.

\begin{Thm}\label{thm:bRdissipative}
Let $W = \indtr W$ be a Gaussian zero--mean process with stationary increments
and continuous paths. If 
\equh\label{eq:limwt}
\lim_{|t|\to\infty}\left(W_t-\sigma_t^2/2\right) = -\infty,\ \mbox{ almost surely}
\eque
then,
\equh\label{eq:intewt}
\int_{-\infty}^{\infty} e^{W_t-\sigma^2_t/2}dt <\infty,\ \mbox{ almost surely},
\eque
where $\sigma_t^2 = \esp W_t^2 = \var(W_t)$.
\end{Thm}
\begin{proof} 
Let $\indtr X$ be the Brown--Resnick process defined
in~\eqref{rep:brownResnick}. Note that the process $\indtr {\log X}$ is also
max--stable but it has Gumbel marginals. Kabluchko {\it et al.}\cite{kabluchko08stationary} have
shown that $\indtr{\log X}$ is stationary and hence so is $\indtr X$. Moreover, by Theorem~13
in~\cite{kabluchko08stationary}, Condition~\eqref{eq:limwt} implies
that $\indtr{\log X}$, or equivalently, $\indtr X$, has a mixed moving
maxima representation. On the other hand, Theorem~\ref{thm:mmm}
implies that any process with mixed moving maxima representation is
dissipative. Dissipativity of $\indtr X$ is equivalent
to~\eqref{eq:intewt} by Theorem~\ref{thm:cons-diss}. This completes
the proof.
\end{proof}

\noindent The following question arises.

\begin{Question}  For what general classes of continuous--path, zero mean
Gaussian processes $\{W_t\}_{t\in\bbR}$ with stationary increments, is the Brown--Resnick 
stationary process \eqref{rep:brownResnick} {\it purely dissipative}?
\end{Question}

\noindent
The next result provides a {\it partial} answer to this question for the interesting case
when $W=\{W_t\}_{t\in\bbR}$ is the  fractional Brownian motion (fBm).  Recall that the fBm is
a zero--mean  Gaussian processes with stationary increments, which is self--similar.  The process $W$ is
said to be self--similar with self--similarity parameter $H>0$, if for all $c>0$, we have that 
$\{W_{ct}\}_{t\in \bbR}\stackrel{{\rm d}}{=}\{c^H W_t\}_{t\in\bbR}$. The fBm necessarily has the
covariance function
\begin{equation}\label{eq:fBm}
 \E W_t W_s = \frac{\sigma^2}{2} {\Big(} |t|^{2H} + |s|^{2H} - |t-s|^{2H}{\Big)},\ \mbox{ with } t,\ s\in\bbR,
\end{equation}
where $0<H\le 1$ is the self--similarity parameter of $W$.
The fractional Brownian motions have versions with continuous paths (see e.g.\ \cite{samorodnitsky94stable}).

\begin{Prop} The stationary Brown--Resnick processes $X=\{X_t\}_{t\in\bbR}$ associated with the
fractional Brownian motions $\{W_t\}_{t\in\bbR}$ in \eqref{eq:fBm} are purely dissipative and
hence they have mixed moving maxima representations.
\end{Prop}
\begin{proof} Without loss of generality, we will suppose that the fBm $W$ has continuous paths.  
As indicated above, the stationarity of $X$ follows from the fact that $W$ has
stationary increments (see Kabluchko et al.~\cite{kabluchko08stationary}).  Now, by 
Theorem~\ref{thm:cons-diss}, $X$ is dissipative, if and only if
\equh\label{eq:ebt}
 \int_{-\infty}^\infty \exp\left\{W_t-\sigma_t^2/2\right\}\d t <\infty,\ \mbox{ almost surely}.
\eque
It is enough to focus on the integral $\int_{0}^\infty \exp\left\{W_t-\sigma_t^2/2\right\}\d t$. By
the Law of the Iterated Logarithm for fractional Brownian motion
(see Oodaira~\cite{oodaira72strassen}), we have
\[
 \limsup_{t\to\infty}W_t/\sqrt{2\sigma_t^2\log\log t} = 1,\ \mbox{ almost surely}.
\]
Hence, with probability one, for any $\delta>0$, there exists $T_1$ (possibly random)
such that $\forall t>T_1$, we have $W_t<(1+\delta)\sqrt{2\sigma^2_t\log\log t}$ almost surely. 
Moreover, there exists $T_2$ sufficiently large (possibly random), such that $\forall t>T_2$, we
have 
$$
 (1+\delta)\sqrt{2\sigma^2_t\log\log t}< \sigma_t^2/4 \equiv \sigma^2 t^{2H}/4\,\malmosts,
$$
where $H \in (0,1]$ is the self--similarity parameter of the fractional Brownian motion $W$.
Now, let $T_0 = \max(T_1,T_2)$. It follows that 
\[
\int_{T_0}^\infty \e^{W_t-\sigma_t^2/2}\d t < \int_{T_0}^\infty \e^{(1+\delta)\sqrt{2\sigma^2_t\log\log t} -
   \sigma^2 t^{2H}/ 2}\d t \leq \int_{T_0}^\infty \e^{-\sigma^2 t^{2H}/4}\d t<\infty\,\malmosts,
\]
which implies \eqref{eq:ebt} since $W_t$ is continuous, with probability one..
\end{proof}

Observe that the above result continues to hold even in the degenerate case $H=1$.  One then
has that 
$W_t = tZ,\ t\in \bbR$, where $Z$ is a zero--mean Gaussian random variable.  In this case, 
the corresponding Brown--Resnick process has a simple {\it moving maxima} representation.  
Indeed, for simplicity,
let $\sigma^2 = {\rm Var}(Z) = 1$ and observe that
$$
 X_t := \Eint{E} e^{tZ(u) - t^2/2} M_1(du) = \Eint{E} e^{Z^2(u)/2} e^{- (t-Z(u))^2/2} M_1(du). 
$$
Note that the measure $\nu(A):= \int_E\ind_{\{Z(u)\in A\}} e^{Z^2(u)/2} \mu(du) \equiv \lambda(A)/\sqrt{2\pi}$, is up to a constant factor equal to the
Lebesgue measure $\lambda$ on $\bbR$.  Therefore, one can show that
$$
\{X_t\}_{t\in \bbR} \eqd {\Big\{} \frac{1}{\sqrt{2\pi}} \Eint{\mathbb R} e^{-(t-z)^2/2} \wtilde M_1(dz){\Big\}}_{t\in \bbR},
$$
where $\wtilde M_1$ is a $1-$Fr\'echet random sup--measure with the Lebesgue control measure.  This shows that $X$ in this simple case 
is merely a {\it moving maxima} rather than a {\it mixed moving maxima}.

We have thus shown that the Brown--Resnick process~\eqref{rep:brownResnick} driven by fractional Brownian motion $\indt W$ is purely dissipative.  
Thus, by Theorem~\ref{thm:mmm} we have that $\indt X$ is a {\it mixed moving maxima}.  It is not clear how one can prove this fact without the use of 
our classification results.  In two very recent papers~\cite{kabluchko08stationary,kabluchko08spectral}, Kabluchko and co--authors established very similar classification results by using very different methods based on Poisson point processes on abstract path--spaces.  Their approach 
yields directly the moving--maxima representation (and hence dissipativity) of the Brown--Resnick type processes $X$ under the alternative 
Condition \eqref{eq:limwt}.  This condition is only shown to be {\it sufficient} for dissipativity of $X$.  Its relationship with
our {\it necessary and sufficient condition} \eqref{eq:intewt} is a question of independent interest.  

The question raised in Kabluchko \cite{kabluchko08spectral} on whether there exist stationary 
Brown--Resnick processes $X$ of mixed type i.e.\ with non--trivial dissipative and conservative
components still remains open.  In view of our new necessary and sufficient condition 
\eqref{eq:intewt}, this question is {\it equivalent} to the following: 

\begin{Question} Is it true for Gaussian processes $W= \{W_t\}_{t\in\bbR}$ with stationary
increments and continuous paths that 
$\mu{\{} \int_{-\infty}^{\infty} e^{W_t-\sigma^2_t/2}dt <\infty {\}} \in \{0,1\}$?

\end{Question}
\appendix
\section{Proofs and Auxiliary Results}
\subsection{Proofs and auxiliary results for Section~\ref{sec:maxLinear}}\label{sec:proofMaxLinear}
The proof of Theorem~\ref{hardin81thm:4.2p}, as well as the auxiliary results in Section~\ref{sec:maxLinear}, follow closely the proofs in Hardin~\cite{hardin81isometries}. There the author dealt with linear isometries instead of max--linear isometries.

\begin{Lem}\label{lem:maxLinearExtension}
Let $(S_1,\calS_1,\mu_1)$ and $(S_2,\calS_2,\mu_2)$ be two measure spaces and suppose that
$T:\calS_1 \to \calS_2$ is a regular set isomorphism.\\
\itemnumber {i} The regular set isomorphism $T$ induces a canonical function map $Tf$, defined mod $\mu_2$ 
for all measurable $f$ on $(S_1,\calS_1)$ and such that $\{Tf \in B\} = T \{f\in B\}$, mod $\mu_2$, 
for all Borel sets $B\in \calB_ {\mathbb R}$. 

The function map $T$ is monotone, linear, max--linear, and preserves the convergence almost everywhere, modulo null sets.
Moreover, $T(f_1 f_2) = T(f_1) T(f_2)$, mod $\mu_2$, for all measurable functions $f_1$ and $f_2$ on $(S_1,\calS_1)$.\\
\itemnumber {ii} If the regular set isomorphism $T$ is measure preserving, then the induced function map $T$
is a max--linear isometry from $\lap(S_1,\calS_1,\mu_1)$ to $\lap(S_2,\calS_2,\mu_2)$, for every $\alpha>0$.
Furthermore, if $T$ is onto, then so is the induced max--linear isometry.
\end{Lem}
\begin{proof}
\itemnumber {i} Consider the sets $A_r :=\{ f\le r\},\ r\in {\mathbb Q}$, where ${\mathbb Q}$ denotes the set of rational numbers.
By the monotonicity of the regular set isomorphism, we obtain that $TA_r \subset TA_s$,  mod $\mu_2$,
for all $r<s,\ r,s\in{\mathbb Q}$. Let $(Tf)(x):= \inf\{ s\, :\, x\in TA_s\}$ and observe that $Tf$ is measurable. Indeed,$$
 \{ Tf \le r\} = \bigcap_{s > r,\ s\in{\mathbb Q}} TA_s = T A_r =  T\{f\le r\},\ \ \ \mbox{ for all }r\in{\mathbb Q}. 
$$
We will show next that $T\{f\in B\} = \{Tf \in B\}$ mod $\mu_2$, for all Borel sets $B \subset {\mathbb R}$.  Indeed, consider
the class of sets:
$$
{\mathcal D} := \{ B \in {\cal B}_{\mathbb R}\, :\, T\{f\in B\} = \{Tf \in B\},\ \mbox{ mod $\mu_2$ }\},
$$
and observe that ${\cal C} \subset {\cal D}$, where ${\mathcal C} := \{ (-\infty,r]\, :\, r\in \mathbb Q\}$.
One can show that ${\cal D}$ is a $\sigma$--algebra and therefore $\sigma({\cal C}) \subset {\cal D}$. This however
implies that ${\cal B}_{\bbR} \equiv {\cal D}$, since ${\cal B}_{\bbR} \equiv \sigma({\cal C})$, thereby showing that
$T\{f\in B\} = \{Tf \in B\},$ for all $B\in {\cal B}_{\bbR}$.

By the properties of regular set isomorphism, it is easy to see that for any sequence
of measurable functions $\{f_m\}_{m\in\bbN}$,
$T(\sup_{m\in\bbN} f_m) = \sup_{m\in\bbN} T f_m$, and $T(\inf_{m\in\bbN} f_m) = \inf_{m\in\bbN} T f_m$, mod $\mu_2$.  Therefore, $T$ preserves pointwise limits of measurable functions, modulo null sets.

One clearly has that $T(\lambda f) = \lambda Tf,$ for all $\lambda\in\bbR$ and measurable $f$'s.  The 
linearity of $T$ follows then from the fact that
$$
 T\{ f+g \le r\} = T {\Big(} \cup_{s\in\bbQ} \{f\le r-s\} \cap \{ g\le s\} {\Big)} = \cup_{s\in\bbQ} \{Tf\le r-s\} \cap 
 \{ Tg\le s\},
$$
which equals $\{Tf +Tg \le r\}.$
The max--linearity of $T$ can be established similarly. The fact that $T$ preserves products, i.e.\ 
$T(f_1 f_2) = T(f_1) T(f_2)$ mod $\mu_2$, for measurable $f_1$ and $f_2$ can be established similarly for non--negative
functions, and then shown to hold for arbitrary functions, by linearity.\\
\itemnumber {ii} Now, if $T$ is measure preserving, then for all simple functions $f = \bigvee_{i=1}^n \lambda_i \ind_{A_i},$ 
with disjoint $A_i$'s in $\calS_1$, we have, by max linearity, that $Tf = \bigvee_{i=1}^n \lambda_i \ind_{TA_i}$. 
Thus, since the regular set isomorphism $T$ is measure preserving, we have
\eqnh
\left\|T(\bigvee\lambda_i\ind_{A_i})\right\|_{\lapimN} = \left\|\bigvee\lambda_i\ind_{T(A_i)}\right\|_{\lapimN} = \left\|\bigvee\lambda_i\ind_{A_i}\right\|_{\lapdoN} 
\,.
\eqne
Now, for any  $f\in\lap(S_1,\calS_1,\mu_1)$, let $\indt f\subset \lap(S_1,\calS_1,\mu_1)$ be a monotone sequence of simple functions such that $f_n\nearrow f\mod\mu_1$ as $n\to\infty$.
We then have that $Tf_n \nearrow Tf$, mod $\mu_2$, as $n\to\infty$, and
 $\int_{S_2}(Tf_n)^\alpha d\mu_2 = \int_{S_1} f_n^\alpha d\mu_1 \nearrow 
\int_{S_1} f^\alpha d\mu_1 <\infty$, as $n\to\infty$.  The monotone convergence theorem implies that $Tf\in \lap(S_2,\calS_2,\mu_2)$ and
 $\|Tf\|_{\lapimN} = \|f\|_{\lapdoN}$, which shows that $T$ is a max--linear isometry.
 If $T$ is {\it onto} as a regular set isomorphism, then the induced max--linear isometry is clearly {\it onto} the set of all
 simple functions, and therefore $T$ is {\it onto} $\lap(S_2,\calS_2,\mu_2)$. \ifthenelse{\boolean{qedTrue}}{\qed}{}
\end{proof}

\noindent
The following lemma is used repeatedly in the sequel. 
\begin{Lem}\label{lem:Rudin}
Let $\filF$ be a max--linear sub--space of $\lapdo$, where $\mu_1$ is a finite measure. Let $U:\filF\to\lapim$ be a max--linear isometry. 
If $\inddo\in\filF$ and $U\inddo = \indim$, then for any collection of functions $\indn f\subset\filF$, 
\equh\label{eq:equimeasurability}
\mu_1(\{f_1,f_2,\dots\}\in B) = \mu_2(\{Uf_1,Uf_2,\dots\}\in B)\,,\quad \forall B\in\calB_{\mathbb R_+^{\mathbb N}}\,.
\eque
Here $\calB_{\mathbb R_+^{\mathbb N}}$ denotes the Borel $\sigma$-algebra on the product space $\mathbb R_+^{\mathbb N} = [0,\infty)^{\mathbb N}$.
\end{Lem}
\begin{proof}
 Since $\indim\in \lapim$, it follows that $\mu_2$ is a finite measure. Without loss of generality, suppose $\mu_1$ and $\mu_2$ are probability measures. 
Let $\mu_{1,\infty}(B) = \mu_1(\{f_1,f_2,\dots\}\in B)$ and $\mu_{2,\infty}(B) = \mu_2(\{Uf_1,Uf_2,\dots\}\in B)$, $\forall B\in\calB_{\mathbb R_+^{\mathbb N}}$ and $f_i\in\filF\,,i\in\mathbb N$. 
In order to show $\mu_{1,\infty} = \mu_{2,\infty}$, we first show that $\mu_{1,\infty}$ and $\mu_{2,\infty}$ induce the same measure $\mu_{1,n}$ and $\mu_{2,n}$ on $\mathbb R_+^n$ via $\mu_{i,n}(B_n) = \mu_{i,\infty}(\widebar B_n)$ for $i=1,2$, where $B_n\in\calB_{\mathbb R_+^n}$ and $\widebar B_n = \{x = (x_1,x_2,\dots)\in\mathbb R_+^{\mathbb N}:(x_1,\dots,x_n)\in B_n\}$. Indeed,
by using a change of variables and the fact that $U$ is a max--linear isometry, we obtain that for all $n\in\mathbb N$, $a_i>0,f_i\in\filF,1\leq i\leq n$,
\begin{multline*}
 \int_{\rnp}\ind \vee \Big(\bigvee_{1\leq i \leq n} a_iz_i\Big)^\alpha d\mu_{1,\infty}(z) = \int_{S_1}\ind \vee \Big(\bigvee_{1\leq i \leq n} a_if_i(x)\Big)^\alpha d\mu_1(x) \\
= \int_{S_2}\ind \vee \Big(\bigvee_{1\leq i \leq n} a_iUf_i(y)\Big)^\alpha d\mu_2(y) =  \int_{\rnp}\ind \vee \Big(\bigvee_{1\leq i \leq n} a_iz_i\Big)^\alpha d\mu_{2,\infty}(z)\,,
\end{multline*}
which implies $\mu_{1,n}(B_n) = \mu_{2,n}(B_n),\forall B_n\in\calB_{\mathbb R_+^n}$. By Lemma 4.1 in~\cite{dehaan86stationary}, $\mu_{1,n} = \mu_{2,n}\,,\forall n\in\mathbb N$. That is, $\mu_{1,\infty}$ and $\mu_{2,\infty}$ agree on the field of all cylinder sets of $\mathbb R_+^\infty$. By the Carath\'eodory extension theorem, $\mu_{1,\infty} = \mu_{2,\infty}$.
\ifthenelse{\boolean{qedTrue}}{\qed}{}
\end{proof}

\begin{proof}[Proof of Theorem~\ref{hardin81thm:2.2p}]
First observe that $\mu_1$ and $\mu_2$ are finite measures, since $\inddo\in\lapdo$ and $\indim\in\lapim$.
We start by defining $T$ on a separable $\sigma$--field and then we verify the consistency of the definition. Let $C = \indn f\subset\filF$ be a countable collection of functions. We then have that 
\equh\label{eq:sigmaC}
\sigma(C) = \left\{\{(\inftydots f)\in B\}\,, B \in \calB_{\mathbb R_+^{\mathbb N}}\right\}\,,
\eque
where $\sigma(C)$ is the minimal $\sigma$-algebra generated by $C$. Observe that for any $A\in\sigma(C)$, there exists $B_A\in\calB_{\mrn_+}$ such that $A = \{(\inftydots f)\in B_A\}$. We therefore define $T_C:\sigma(C)\to\sigma(U(C))$ as
\equh\label{eq:TC}
T_CA = \left\{(\inftydots {Uf})\in B_A\right\}\,.
\eque
In the sequel, we will first show that $T_C$ is: (1) well defined (note that $B_A$ may not be unique for a given $A$), (2) measure preserving, (3) onto and (4) $T_C$ induces a unique max--linear isometry from $\lapdoSC$ onto $\lapimSC$, implying that $T_C$ satisfies (i) and (ii). Then we show that for the induced max--linear isometry $T_C$: (5) $T_C$ and $U$ coincides on $\cmsspan(f_t,t\in T)$ and (6) $T_C$ is unique, implying that $T_C$ satisfies (iii).\\
\itemnumber1 \textit{$T_C$ is well defined modulo $\mu_2$--null sets.} Indeed, if there is another $\widetilde B_A\in\calB_{\mrn_+}$ such that $A = \{(\inftydots {f})\in \widetilde B_A\}$, then
\begin{multline*}
\mu_2\Big(\{(\inftydots {Uf})\in B_A\}\Delta\{(\inftydots {Uf})\in \widetilde B_A\}\Big) = \mu_2\Big(\{(\inftydots {Uf})\in B_A \Delta \widetilde B_A\}\Big)\,.
\end{multline*}
By Lemma~\ref{lem:Rudin}, the last expression equals
\eqnh
& & \mu_1\Big(\{(\inftydots {f})\in B_A \Delta \widetilde B_A\}\Big) \\
\eqnhspace = \mu_1\Big(\{(\inftydots {f})\in B_A\}\Delta\{(\inftydots {f})\in \widetilde B_A\}\Big) = \mu_1(A\Delta A) = 0\,,
\eqne
which shows that the definition of $T_CA$ in~\eqref{eq:TC} does not depend on the choice of $B_A\in\calB_{\mrn_+}$, modulo $\mu_2$--null sets.\\
\itemnumber2 \textit{$T_C$ is a regular set isomorphism.} For convenience, we prove the three conditions in Definition~\ref{def:regular} in different order. For every $A\in\sigma(C)$, let $B_A \in \calB_{\mrn_+}$ be such that $\{(\inftydots f)\in B_A\} = A$. First, by Lemma~\ref{lem:Rudin}, we have
\equh\label{eq:measurePreserving}
\mu_1(A) = \mu_1(\{(\inftydots f)\in B_A\} )= \mu_2(\{(\inftydots {Uf}) \in B_A\}) = \mu_2(T_C(A))\,, 
\eque
which shows that (iii) of Definition~\ref{def:regular} holds.
Relation~\eqref{eq:measurePreserving} shows moreover that $T_C$ is \textit{measure-preserving}. Second, note that $S_1\in\sigma(C)$. Then we can show (i) of Definition~\ref{def:regular} by 
\eqnh
& & T_C(S_1\setminus A) = \{(\inftydots {Uf})\in \mrn_+\setminus B_A\} \\
& & \ \ \ \ = \{(\inftydots {Uf})\in \mrn_+\}\setminus\{(\inftydots {Uf}) \in B_A\} = T_C(S_1)\setminus T_C(A)\,. 
\eqne
Finally, suppose $\inftydots A$ are arbitrary disjoint sets in $\sigma(C)$. Observe that, modulo $\mu_1$, 
\[
\cup_{n=1}^\infty A_n = \cup_{n=1}^\infty\left\{\right(\inftydots f)\in B_{A_n}\}
= \left\{\right(\inftydots f)\in \cup_{n=1}^\infty B_{A_n}\}\,,
\]
whence, $\cup_{n=1}^\infty B_{A_n}$ may be viewed as a particular choice of $B_{\cup_{n=1}^\infty A_n}$, and thus, in view of~\eqref{eq:TC},
\eqnh
T_C(\cup_{n=1}^\infty A_n)
 & = & \{(\inftydots {Uf}) \in \cup_{n=1}^\infty B_{A_n}\}
 \\
 & = & \cup_{n=1}^\infty \{(\inftydots {Uf}) \in B_{A_n}\}
 = \bigcup_{n=1}^\infty\{T_C(A_n)\}\,,
\eqne
implying (ii) of Definition~\ref{def:regular}.

\itemnumber3 \textit{$T_C$ is onto $\sigma(U(C))$.} Observe first that as in~\eqref{eq:sigmaC}, $\sigma(U(C)) = \{(\inftydots {Uf})\in \calB_{\mrn_+}\}$. Now $\forall U_A\in\sigma(U(C))$ there is a $B_A \in\calB_{\mrn_+}$ such that $U_A = \{(\inftydots {Uf})\in B_A\}$.
Then, let $A=\{(\inftydots f)\in B_A\}$ and observe that by~\eqref{eq:TC}, we have $T_CA = \{(\inftydots{Uf})\in B_A\}=U_A$. This shows that $T_C$ is onto.\\
\itemnumber4 \textit{$T_C$ induces a max--linear isometry from $\lapdoSC$ onto $\lapimSC$.} The proof is standard and is given as Lemma~\ref{lem:maxLinearExtension} in Appendix~\ref{sec:proofMaxLinear}.
\\
\itemnumber5 \textit{$T_C$ and $U$ coincide on $\cmsspan\{\inddo,\inftydots f\}$.} Observe that $T_C\inddo = \indim$ and for any $B\in\calB_{\mathbb R_+}$, 
\[
\{T_Cf_j\in B\} = T_C(\{f_j\in B\}) = \{Uf_j\in B\}\,, j=1,2,\dots\,,
\]
where the first equality follows as in the case of linear mapping (e.g. see p452-454~\cite{doob53stochastic}) and the second equality follows by~\eqref{eq:TC}.
This shows that $T_C(f_j) = U(f_j)\,,j=1,2,\dots$. Since $T_C$ is a max--linear isometry by (4), $T_C$ and $U$ coincide on any finite positive max--linear combinations of $\inddo,\inftydots f$ and hence on $\cmsspan\{\inddo,\inftydots f\}$.\\
\itemnumber6 \textit{$T_C$ is unique.} Assume that there exists another max--linear isometry $V$ from $\lapdoSC$ to $\lapim$ such that $V$ and $U$ agree on $\cmsspan\{\inddo,\inftydots f\}$. 
We will show that $V$ and $T_C$ coincide on $\lapdoSC$. It is enough to show that, for any $A_0\in\calB_{\mathbb R_+}$ and any $A = \{(\inftydots {Uf})\in\calB_A\}\in\sigma(U(C))$, we have
\[
\mu_2(\{T_Cf\in A_0\}\cap A) = \mu_2(\{Vf\in A_0\}\cap A)\,.
\]
Indeed, for any $f\in\lapdoSC$ we have, by (5),  
\eqnhn
& & \mu_2(\{T_Cf\in A_0\}\cap A) = \mu_2(\{(T_Cf,\inftydots {T_Cf})\in A_0\times B_A\}) \nonumber\\
\eqnhspace = \mu_1(\{(f,\inftydots f)\in A_0\times B_A\})\label{eq:nuTC}\\
\eqnhspace = \mu_2(\{(Vf,\inftydots {Vf})\in A_0\times B_A\}) \nonumber\\
\eqnhspace = \mu_2(\{(Vf,\inftydots {Uf})\in A_0\times B_A\})\label{eq:nuVf}\\
\eqnhspace  = \mu_2(\{Vf\in A_0\}\cap A)\nonumber
\eqnen
by using Lemma~\ref{lem:Rudin} in~\eqref{eq:nuTC} and~\eqref{eq:nuVf}. Hence, we have $\{T_Cf\in A_0\} = \{Vf\in A_0\}\,,\forall A_0\in\calB_{\mathbb R_+}$ and $\forall f\in\lapdoSC$, which implies that $T_C$ and $V$ agree on $\lapdoSC$.

To complete the proof, we define a max--linear isometry from
$\lapdoS$ {\it onto} $\lapimS$.  Note that, for all
$f\in\lap(S_1,\sigma(\filF),\mu_1)$ there exists a countable
collection of functions $C = \{\inftydots f\}$, such that
$f\in\sigma(C)$. We therefore define $Tf = T_Cf$.  To check the
consistency of this definition, suppose that $f\in\sigma(\widetilde
C)$, for another countable collection of functions $\widetilde
C\subset \filF$.  Since $C\subset C\cup\widetilde C$, by
using~\eqref{eq:TC} one can show that $T_C(A) = T_{C\cup\widetilde
C}(A)$ for every $A\in\sigma(C)$. Thus, $T_Cf = T_{C\cup\widetilde
C}f$ and similarly $T_{\widetilde C}f = T_{C\cup\widetilde C}f$, which
shows that $T$ is well-defined.  It is easy to see that $T$ is induced
by a measure preserving regular set isomorphism of $\sigma(\filF)$
\textit{onto} $\sigma(U(\filF))$. This is because that for every
$A\in\sigma(\filF)$, we have $\ind_A\in\sigma(C)$ with some countable
collection $C\subset\filF$.
\end{proof}

\begin{proof}[Proof of Lemma~\ref{lem:ratio1}]
First, to show $\rho(F) = \rho(\cmsspan(F))$, it suffices to show that $\rho(\cmsspan(F))\subset\rho(F)$. Observe that 
for any $f_i,g_i\in F,a_i\geq 0,b_i\geq 0,i\in\mathbb N$ and $c>0$, we have
\begin{multline*}
\Big\{\frac{\bigvee_{i\in\mathbb N}a_if_i}{\bigvee_{j\in\mathbb N}b_jg_j}\leq c\Big\} = \bigcap_{i\in\mathbb N}\Big\{\frac{a_if_i}{\bigvee_{j\in\mathbb N} b_jg_j}\leq c\Big\}\\ = \bigcap_{i\in\mathbb N}\bigcap_{n\in\mathbb N}\Big\{\frac{a_if_i}{\bigvee_{j\in\mathbb N} b_jg_j}<c+\frac1n\Big\} = 
\bigcap_{i\in\mathbb N}\bigcap_{n\in\mathbb N}\bigcup_{j\in\mathbb N}\Big\{\frac{a_if_i}{b_jg_j}< c+\frac1n\Big\} \,.
\end{multline*}
Hence $\rho(\cmsspan(F))\subset\rho(F)$. Next, $\rho(F)\subset\sigma(F)$ follows from the fact that 
\[
\left\{f_1/f_2\leq c\right\} = \bigcup_{0 < q\in\mathbb Q}\left(\left\{f_1 \leq q\right\}\cap\left\{f_2\geq q/c\right\}\right)\,,
\]
for any $f_1,f_2\in F$ and $c>0$, where $\mathbb Q$ denotes the set of rational numbers. Finally, when $\ind_S\in F$, $\rho(F)\supset\rho(f/\ind_S:f\in F) = \sigma(F)$. As it is always true that $\rho(F)\subset\sigma(F)$, it follows $\rho(F) = \sigma(F)$.
\end{proof}

\begin{proof}[Proof of Lemma~\ref{lem:fullSupport1}]
Suppose first that $\filF$ is separable, i.e., there exists a countable collection of functions $f_n\in\filF,n\in\mathbb N$, such that $\filF = \cmsspan\{\inftydots f\}$. Then, we have that
\[
g = \bigvee_{n=1}^\infty \left(f_n/({2^n\left\|f_n\right\|_{\lap(S,\mu)}})\right)\mbox{ belongs to }\laps\,,
\]
because
\[
\left\|g\right\|_{\lap(\mu)}^\alpha = \int\bigvee_{n=1}^\infty\frac{f_n^\alpha}{2^{n\alpha}\left\|f_n\right\|_{\laps}^\alpha}d\mu\leq \sum_{n=1}^\infty \int\frac{f_n^\alpha}{2^{n\alpha}\left\|f_n\right\|_{\laps}^\alpha} d\mu = \sum_n\frac1{2^{n\alpha}} <\infty\,.
\]
Since $\filF$ is a max--linear space, we have $g\in\filF$ and clearly $g$ has full support in $\filF$ since for any $f\in\filF$, $\supp(f)\subset\bigcup_{n=1}^\infty\supp(f_n) = \supp(g)\mod\mu$. 

Next consider the case when $\mu$ is $\sigma$-finite. Let $\widebar\mu$ be a finite measure equivalent to $\mu$ (i.e., $\widebar\mu\ll\mu$ and $\mu\ll\widebar\mu$). Now let $F\subset\filF$ be any arbitrary countable collection of functions in $\filF$, set $s(F)\defe \widebar\mu\left(\bigcup_{f\in F}\supp(f)\right)$ and define $s\defe \sup_{F\in\filF} s(F)$. 
Thus, consider a sequence $F_n\subset\filF,n\in\mathbb N$ of countable collections of functions, such that $s(F_n)\uparrow s$ as $n\to\infty$. Let $C = \bigcup_{n\in\mathbb N} F_n$ and observe that $C$ is countable. Then by the first part of the proof, there exists $g\in\cmsspan(C)$ with full support in $\cmsspan(C)$, since $\widebar\mu\left(\bigcup_{f\in C}\supp(f)\setminus\supp(g)\right) = 0$ implies $\mu\left(\bigcup_{f\in C}\supp(f)\setminus\supp(g)\right) = 0$. The function $g$ has also full support in $\filF$. Indeed, if there exists a function $f_0\in\filF$ such that $\widebar\mu\left(\supp(f_0)\setminus\supp(g)\right) =\epsilon > 0$, then $f_0\notin C$ and $\lim_{n\to\infty}s(F_n\cup\{f_0\}) \geq s+\epsilon > s$, which is a contradiction. This completes the proof of the lemma.\ifthenelse{\boolean{qedTrue}}{\qed}{}
\end{proof}

\begin{proof}[Proof of Lemma~\ref{lem:fullSupport2}]
Let $g_0 = Uf_0$ and let $g_1 = Uf_1$ for an arbitrary $f_1\in\filF$. We clearly have that $f_2\defe f_0\vee f_1$ and $g_2\defe g_0\vee g_1 = Uf_2$ have full supports in $\cmsspan\{f_1,f_2\}$ and $\cmsspan\{g_1,g_2\}$, respectively. To prove the result, it is enough to show that $\mu_2(\supp(g_1)\setminus \supp(g_0)) = 0$, or equivalently, $\mu_2(\supp(g_2)\setminus\supp(g_0)) = 0$.

Consider the finite measures
\equh\label{eq:nu1}
\nu_1 = f_2^\alpha d\mu_1\quad\mbox{and}\quad \nu_2 = g_2^\alpha d\mu_2\,,
\eque
restricted to the spaces $(\supp(f_2),\left.\calB_{S_1}\right|_{\supp(f_2)})$ and $(\supp(g_2),\left.\calB_{S_2}\right|_{\supp(g_2)})$. Now, define
\[
V\left(a\ind_{\supp(f_2)}\vee b(f_0/f_2)\right)\defe a\ind_{\supp(g_2)}\vee b(g_0/g_2)\quad \forall a,b\geq 0\,.
\]
Observe that
\eqnh
& & \int_{\supp(f_2)}\Big(a\ind\vee\lambda bf_0\frac1{f_2}\Big)^\alpha d\nu_1 = 
\int_{S_1}\Big(af_2\vee\lambda bf_0\Big)^\alpha d\mu_1 \\
\eqnhspace \ \ = \int_{S_2}\Big(ag_2\vee\lambda bg_0\Big)^\alpha d\mu_2 = 
\int_{\supp(g_2)}\Big(a\ind\vee\lambda bg_0\frac1{g_2}\Big)^\alpha d\nu_2\,.
\eqne
This shows that $V:\cmsspan\{\ind_{\supp(f_2)},f_0/f_2\}\to\lap(\supp(g_2),\nu_2)$ is a max--linear isometry mapping $\ind_{\supp(f_2)}$ to $\ind_{\supp(g_2)}$. Thus, by Lemma~\ref{lem:Rudin}, we obtain
\[
\nu_1\left(\{f_0/f_2 = 0\}\right) = \nu_2\left(\{V(f_0/f_2)= 0\} \right) = \nu_2\left(\{g_0/g_2 = 0\}\right) = 0\,,
\]
Since $\nu_1(\{f_0/f_2=0\}) = \nu_1(\supp(f_2)\setminus\supp(f_0)) = 0$, we also have that $\nu_2(g_0/g_2=0) = \nu_2(\supp(g_2)\setminus\supp(g_0)) = 0$. This, in view of~\eqref{eq:nu1}, implies that $\mu_2(\supp(g_2)\setminus\supp(g_0)) = 0$ and hence $\mu_2(\supp(g_1)\setminus\supp(g_0)) = 0$, since $\supp(g_1)\subset\supp(g_2)$. We have thus shown that for an arbitrary $f_1\in\filF$, $\mu_2(\supp(Uf_1)\setminus\supp(g_0)) = 0$, which shows that $Uf_0 = g_0$ has full support in $U(\filF)$ (see Definition~\ref{def:fullSupport}).
\ifthenelse{\boolean{qedTrue}}{\qed}{}
\end{proof}
\begin{proof}[Proof of Theorem \ref{hardin81thm:4.2p}]
Let $f_0\in\calF$ be a function with full support in $\filF$, i.e., $\supp(f_0) = S_1$ 
(Lemma~\ref{lem:fullSupport1}). Define $\filF_0\defe \left\{f\cdot(1/f_0)\,,f\in\filF\right\}$. 
 Since $f_0/f_0=\inddo\in\filF_0$, it follows that 
\begin{equation}\label{eq:sig-f0}
 \sigma(\filF_0) = \rho(\filF_0) = \rho(\filF),
\end{equation}
where the second equality follows from the fact that $f_1/f_2 = (f_1/f_0)/ (f_2/f_0)$ mod $\mu_1$, for all
$f_1, f_2 \in \calF$, since $f_0$ has full support in $\calF$.  Therefore, any element $rf\in \eratiosf$, 
$r\in\calR_+(\filF)$ and $f\in\filF$, can be represented as follows:
$$
 rf = \left(rf\cdot(1/{f_0})\right)f_0 = r_0f_0,
$$
where $r_0 = rf\cdot(1/{f_0}) = r\cdot(f/f_0)$ is a $\rho(\filF)$-measurable and hence $\sigma(\filF_0)$-measurable
function. Hence, we have that
\[
\eratiosf = \left\{r_0f_0\in\lapdo\,,r_0\geq 0, r_0 \in\sigma(\filF_0)\right\}\,.
\]

Next, introduce the measures $d\mu_{1,f_0} = f_0^\alpha d\mu_1$ and $d\mu_{2,f_0} = (Uf_0)^\alpha d\mu_2$, and
observe that both of them are finite. We thus have that $\filF_0$ is a max--linear sub--space of
$\lap(S_1,\mu_{1,f_0})$ and similarly $\calG_0 \defe \left\{Uf\cdot(1/{Uf_0})\,, f\in\filF\right\}$ is a
max--linear sub--space $\subset\lap(S_2, \mu_{2,f_0})$.  
It is easy to check that 
\eqnh
U_0:\filF_0 & \to & \calG_0,\ \ \mbox{ defined by } U_0(f) := U(f\cdot f_0)\cdot (1/Uf_0),\ \ f\in{\calF_0}
\eqne
is a max--linear isometry from $\filF_0 \subset \lap(S_1,\mu_{1,f_0})$ to $\calG_0\subset \lap(S_2,\mu_{2,f_0})$.
Note, however, that these two $\lap-$spaces involve finite measures and $U_0\inddo = \ind_{\supp(Uf_0)}$. 
Thus, by Theorem~\ref{hardin81thm:2.2p}, we obtain that $U_0$ has a unique extension to a max--linear isometry
\[
  T:\lap(S_1,\sigma(\filF_0),\mu_{1,f_0})\to\lap(\supp(Uf_0),\sigma(\calG_0)_{\vert_{\supp(Uf_0)}},\mu_{2,f_0})\,,
\]
which is induced by a measure preserving regular set isomorphism $T$ from $\sigma(\filF_0)$ {\it onto} 
$\sigma(\calG_0)$.

We can now construct the desired extension $\widebar U$ of the max--linear isometry $U$.  Consider the mappings
\[
M:\eratiosf\to\lap(S_1,\sigma(\filF_0),\mu_{1,f_0})
\]
defined by $Mf := f\cdot(1/f_0)\,,\forall f\in \eratiosf$ and
$$
N:\lap(\supp(U f_0) ,\sigma(\calG_0)_{\vert_{\supp(U f_0)}},\mu_{2,f_0}) \to\lapim
$$
defined by $Ng := g\cdot (Uf_0), \forall g\in 
\lap(\supp(Uf_0) ,\sigma(\calG_0)_{\vert_{\supp(U f_0)}},\mu_{2,f_0})$.
Note that both mappings $M$ and $N$ are {\it one-to-one} and that $M$ is trivially {\it onto}.  We will now
show that $N$ is also {\it onto}.  Indeed, as in \eqref{eq:sig-f0}, we have that
\begin{equation}\label{e:sig-g0}
 \sigma(\calG_0) = \rho(\calG_0) = \rho(U(\calF)).
\end{equation}
Consider an arbitrary $g \in \eratiosuf,$ and note that $g = r U(f)$, with some $r\in \rho(U(\cal F))$ and
$f\in\calF$.  We have that $g = \widetilde r U(f_0)$ with $\widetilde r = r U(f)/U(f_0)$, since $Uf_0$ 
has full support in $U(\calF)$ (Lemma~\ref{lem:fullSupport2}).  By \eqref{e:sig-g0}, we have that $\widetilde r$
is $\rho(U(\calF))$ and hence $\sigma(\calG_0)-$measurable, and since $g = rU(f) \in\lap(S_2,\mu_2)$, it follows that
$\widetilde r \in \lap(S_2,\sigma(\calG_0),\mu_{2,f_0})$.  This shows that $N(\widetilde r) = r U(f) = g$,
and since $g\in \eratiosuf$ was arbitrary, it follows that $N$ is {\it onto} $\eratiosuf$.

At last, we define
\[
 \widebar U \defe NTM:\eratiosf\to\lapim.
\]
We will complete the proof by verifying that $\widebar U$ satisfies \eqref{hardin81thm:4.2peq:1} and
\eqref{hardin81thm:4.2peq:2} as well as the fact that $\widebar U$ is {\it onto} and {\it unique}. 
To prove \eqref{hardin81thm:4.2peq:1}, observe that 
$$
 \widebar U(rf) =  NTM( rf \cdot (1/f_0)\cdot f_0) = N T(r f \cdot(1/f_0)) =  (Uf_0) T(r) T(f\cdot(1/f_0)),
$$
where the last equality follows from the fact that $T(f_1 f_2) = T(f_1) T(f_2)$, 
for any two measurable functions $f_1$ and $f_2$ (Lemma~\ref{lem:maxLinearExtension}).  Since
 $T(f\cdot (1/f_0)) = U_0(f\cdot(1/f_0)) = U(f)/U(f_0)$, we obtain that
$$
\widebar U(rf) = (U f_0) T(r) U(f)/U(f_0) = T(r) U(f),
$$
which yields \eqref{hardin81thm:4.2peq:1}.

To prove \eqref{hardin81thm:4.2peq:2}, note that for all $A\in \rho(U(\filF))$, we have
\eqnh
\left(\mu_{1,f}\circ T\inv\right) A & = & \int_{S_1} \left(\ind_{T\inv A}f\right)^\alpha d\mu_1 \\
& = &   \int_{S_2} T(\ind_{T^{-1}A})^\alpha T(f)^\alpha d\mu_2  = \int_{S_2} \ind_{A} U(f)^\alpha d\mu_2,
\eqne
which is equivalent to Relation \eqref{hardin81thm:4.2peq:2}.

Now, the extension $\widebar U = NTM$ is {\it onto} $\eratiosuf$ because so are the mappings $M, N$ and $T$.
Finally, to prove the uniqueness of $\widebar U$, suppose that there exists another max--linear isometry, 
$V:\eratiosf\to\lapim$, extending $U$. By the definitions of $M$ and $N$, we have that 
$N\inv V M\inv$ is a max--linear isometry from $\lap(S_1,\sigma(\filF_0),\mu_{1,f_0})$ to
$\lap(\supp(U f_0), \sigma(\calG_0)_{\vert_{\supp(U f_0)}},\mu_{2,f_0})$.  We also have that
\[
  N\inv V M\inv(f/f_0) = N\inv V f = N\inv(Uf) = Uf/(Uf_0)\,,\ \ \ \mbox{ for all }f\in \calF,
\]
which shows that $N\inv V M\inv$ coincides with $U_0$ on $\calF_0$.  Since $U_0$ has a unique extension 
$T:\lap(S_1,\sigma(\calF_0),\mu_{1,f_0})\to \lap(\supp(U f_0),\sigma(\calG_0)_{\vert_{\supp(U f_0)}},\mu_{2,f_0})$,
we obtain that $N\inv VM\inv = T$, which implies $V = NTM \equiv \overline U$.  This completes the proof
of the theorem. \ifthenelse{\boolean{qedTrue}}{\qed}{}
\end{proof}

\begin{proof}[Proof of Lemma~\ref{lem:ratio2}]
Fix $f_0\in F$ with full support. Since $f_0$ is $\sigma(F)$-measurable, so is $1/f_0$. Now for any $f\in \laps$, $f$ is $\sigma(F)$-measurable. Observe that $\rho(F)\subset\sigma(F)\subset\calB_S$, whence, by $\rho(F)\sim\calB_S\mod\mu$, $\rho(F) \sim \sigma(F) \sim \calB_S\mod\mu$. Hence, $f\cdot(1/f_0)$ is $\calB_S$-measurable. Thus $f=(f\cdot(1/f_0))f_0\in\laps$.
\end{proof}

\subsection{Proofs for Sections \ref{sec:classification} and \ref{sec:stationary}}
\label{sec:proofStationary}

\begin{proof}[Proof of Proposition \ref{p:measurability}]  
To prove part {\it (i)}, observe that since $\mu$ is $\sigma$--finite, it is enough to focus
on the case when $\mu$ is a probability measure: $\mu(S)=1$.  Thus, $\{f_t(s)\}_{t\in T}$ may be viewed as a stochastic
process, defined on the probability space $(S,{\cal B}_S,\mu)$.

Note that $L_+^\alpha(S,\mu)$ equipped with the metric $\mrho(f,g)= \int_S |f^\alpha - g^\alpha| d\mu$, is a complete
separable metric space.  Furthermore, $\mrho$ metrizes the convergence in probability in the space $(S,\mu)$.  
Therefore, Theorem 3 of Cohn \cite{cohn72measurable} (see also Proposition 9.4.4 in \cite{samorodnitsky94stable}) implies
that the stochastic process $f=\{f_t(s)\}_{t\in T}$ has a measurable modification {\it if and only if} the map $h_f:t\mapsto [f_t]$ 
is Borel--measurable and has separable range $h_f(T)$. Here $[f]$ denotes the class of all $L_+^\alpha(\mu)$--functions, equal to $f$, $\mu\ae.$ 

Similarly, $X=\{X_t\}_{t\in T}$ has a measurable modification {\it if and only if} $h_X:t\mapsto [X_t]$ is Borel--measurable and has
separable range $h_X(T)$, where $[X_t] \in {\cal L}^0(\Omega,{\cal F},\P)$ is equipped with a metric, which metrizes the
convergence in probability.  Here ${\cal L}^0(\Omega,{\cal F},\P)$ denotes the collection of equivalence classes of random variables,
with respect to the relation of almost sure equality.   We focus on the set ${\cal M} = \{[\xi]\, :\, 
\xi =\eint{S}g d M_\alpha,\ g\in L_+^\alpha(S,\mu)\}$, which is a closed subset of ${\cal L}^0(\Omega,{\cal F},\P)$ with respect to the
convergence in probability.  Theorem 2.1 of \cite{stoev06extremal}, shows that since $(L_+^\alpha(S,\mu),\rho)$ is complete and separable, so
is ${\cal M}$ with respect to the metric:
$$
 \rho_{\cal M}(\xi,\eta) := 2\|\xi\vee \eta\|_\alpha^\alpha - \|\xi\|_\alpha^\alpha - \|\eta\|_\alpha^\alpha.
$$
Furthermore, $\rho_{\cal M}$ metrizes the convergence in probability and we have
\begin{equation}\label{e:rho-equiv}
 \rho_{\cal M}(\xi,\eta) = \int_S |f^\alpha - g^\alpha| d\mu \equiv \rho(f,g),
\end{equation}
for all $\xi = \eint{S} f dM_\alpha$ and $\eta=\eint{S} g dM_\alpha$, with $f, g\in L_+^\alpha(S,\mu)$.

Now, the separability of $L_+^\alpha(S,\mu)$ and ${\cal M}$ implies the separability of 
the ranges $h_f(T)\subset L_+^\alpha(S,\mu)$ and $h_X(T)\subset {\cal M}$, respectively. On the other hand,
the equivalence \eqref{e:rho-equiv} of the two metrics $\rho_{\cal M}$ and $\rho$ implies that $h_f:T\to L_+^\alpha(S,\mu)$ is
Borel--measurable {\it if and only if} $h_X:T\to {\cal M}$ is Borel--measurable.  This, in view of Theorem 3 of Cohn \cite{cohn72measurable}, 
yields {\it (i)}.

In view of Proposition~\ref{prop:conditionS}, to establish {\it (ii)}, we should show that any measurable $\alpha$--Fr\'echet process $X$ satisfies Condition S.
As argued above, the map $h_X:t\mapsto [X_t]$ has a separable range in the metric space ${\cal L}_0(\Omega,{\cal F},\P)$. Hence, there exists a
countable set $T_0\subset T$, such that for all $t\in T$, for some $t_n\in T_0$, we have $X_{t_n}\stackrel{P}{\to} X_t$, as $n\to\infty$.
This shows that the process $X$ is separable in probability (satisfies Condition S, see Definition \ref{d:Cond-S}) and the proof
is complete.
\end{proof}

\begin{proof}[Proof of Theorem \ref{thm:cospectralDecomp}]
Part {\it (ii)} follows immediately from~\eqref{e:S-co-spec}. 
To prove {\it (i)}, consider another measurable representation $\indt{f\topp 2}\subset\lap(S_2,\mu_2)$ of the same
 process $\indt X$. We show that $\indt{f\topp 2}$ also admits a co--spectral decomposition and, letting the corresponding decomposition of the process be
\equh\label{eq:widetildeX}
\indt X \eqd \bccbb{\widehat X\topp 1_t\vee\cdots\vee\widehat X\topp n_t}_{t\in T}\,,
\eque
we have
\equh\label{eq:widetildeX2}
\indt{X\topp j} \eqd \indt{\widehat X\topp j}\,,1\leq j\leq n\,.
\eque
Let $\indt{f\topp 1}\subset\lap(S_1,\mu_1)$ denote the representation in assumption, which admits a co--spectral decomposition w.r.t.\ $\{\calP_j\}_{1\leq j\leq n}$. Without specification, the following arguments hold for both $i = 1,2$. 

First, by Proposition \ref{p:measurability}, the process $X$ has the representation in 
 \eqref{rep:extremalRep}, and hence it has a minimal representation with standardized support
 $\{f_t(s)\}_{t\in T}\subset \lap(S_{I,N},\lambda_{I,N})$ by Theorem \ref{thm:standardized}.  This
 representation can be also chosen to be jointly measurable. By \eqref{eq:pointRep1} in Theorem~\ref{thm:relation}, we have
\equh\label{eq:gts}
f_t^{(i)} (s) = h_i(s)f_t(\Phi_i(s)) =: \widetilde f_t\topp i(s)\,,\mu_i\ae\,,\forall t\in T\,,
\eque
where $h_i:S_i\to\mathbb R_+\setminus\{0\}$ and $\Phi_i$ from $S_i$ onto $\ssS$ are both measurable. 
Since $(t,s)\mapsto f_t(s)$ is measurable,
it follows that $\widetilde f_t^{(i)}(s)$ is jointly measurable modification of $f_t^{(i)}(s)$. Consider the sets
\[
 N^{(i)} \defe \left\{(t,s):f_t^{(i)}(s) \neq \widetilde f_t^{(i)} (s)\right\}\,.
\]
By~\eqref{eq:gts}, we have that $\mu_i(N_t^{(i)}) = 0\,,\forall t\in T$, where 
$N_t^{(i)} = \left\{s:(t,s)\in N^{(i)}\right\}$. Thus, by Fubini's Theorem,
there exists $\widetilde S_i\subset S_i$
such that $\mu_i(S_i\setminus\widetilde S_i) = 0$ and for all
$s\in\widetilde S_i$, $\widetilde f_\cdot^{(i)}(s) =
f_\cdot^{(i)}(s)\,,\lambda\ae$. 

The argument above implies that
\equh\label{eq:ftis}
f_t\topp i(s) = h_i(s)f_t\circ\Phi_i(s), \forall (t,s)\in T\times \widetilde S_i.
\eque

Now, suppose $S_1$ has a co--spectral decomposition $S_1 = \bigcup_{j=1}^nS_1\topp j\mod\mu_1$. We show that this induces a co--spectral decomposition of $\ssS$. Without loss of generality, assume that $S_1\topp j\subset \widetilde S_1, 1\leq j\leq n$. Set
\equh\label{eq:SINj}
\ssS\topp j \defe\Phi_1(S_1\topp j), 1\leq j\leq n\mand \ssS\topp 0 \defe \ssS\setminus\bigcup_{j=1}^n\ssS\topp j\,.
\eque
By~\eqref{eq:ftis}, $\ssS\topp j\subset\{s:f_\cdot(s)\in\calP_j\}, 1\leq j\leq n$. 
Note that the assumption $S_1\topp j\cap S_1\topp k\subset\{s\in S_1: f\topp1_\cdot (s) \equiv 0\}$ implies that $\ssS\topp j\cap\ssS\topp k\subset\{s\in \ssS:f_{\cdot}(s) \equiv 0\}$, for all $1\leq j<k\leq n$. Moreover, $\Phi_1\inv(\ssS\topp 0)\subset S_1\setminus\bigcup_{j=1}^nS_1\topp j$, whence $\ssL(\ssS\topp 0) = 0$. We have thus shown that $\{\ssS\topp j\}_{1\leq j\leq n}$ is a co--spectral decomposition of $\indt f\subset\lap(\ssS,\ssL)$, w.r.t.\ $\{\calP_j\}_{1\leq j\leq n}$.

Next, we show that for any spectral representation $\indt{f\topp 2}\subset\lap(S_2,\mu_2)$, there exists a co--spectral decomposition of $S_2$ w.r.t.\ $\{\calP_j\}_{1\leq j\leq n}$. Indeed, the decomposition is induced by setting $S_2\topp j\defe\Phi_2\inv(\ssS\topp j)\cap\widetilde S_2, 1\leq j\leq n$.
One can easily verify that $\{S_2\topp j\}_{1\leq j\leq n}$ is a co--spectral decomposition w.r.t. $\{\calP_j\}_{1\leq j\leq n}$. 

Finally, by the construction of $\{S_i\topp j\}_{1\leq j\leq n}, i = 1,2$ above, we have
\equh\label{eq:Phii}
\ssL\left(\Phi_i(S_i\topp j)\triangle\ssS\topp j\right) = 0\,, \ \ \forall 1\leq j\leq n\,.
\eque
Note that~\eqref{eq:ftis} induces a max--linear isometry from $\lap(\ssS,\ssL)$ to $\lap(S_i,\mu_i)$. Combining with ~\eqref{eq:Phii} and Remark~\ref{rem:Phi-non-sing}, we have
\[
\Big\{\Eintt_{S_i\topp j}f_t\topp i dM_\alpha\topp i\Big\}_{t\in T} \eqd 
 \Big\{\Eintt_{\ssS\topp j}f_t dM_\alpha\Big\}_{t\in T}\,,1\leq j\leq n\,.
\]
This implies~\eqref{eq:widetildeX2}.
\end{proof}

\begin{proof}[Proof of Theorem~\ref{thm:flow}]
This result can be established by following closely the proof of Theorem~3.1 in~\cite{rosinski95structure} and replacing the linear combination $g_n = \sumin c_{ni}f_{ni}$ therein by the max--linear combination $g_n = \bigvee_{i=1}^n c_{ni}f_{ni}\in\msspan\{f_t:t\in T\}$.
For the completeness, we provide the details next.

Suppose $\indt f$ is minimal. Then, for any $\tau\in T$, by stationarity $\{f_{t+\tau}\}_{t\in T}$ is also a minimal representation of the same $\alpha$--Fr\'echet process. By applying Corollary~\ref{coro:uniquePointMapping}, there exist a one-to-one and onto measurable function $\Phi_\tau:\ssS\to\ssS$ and a measurable function $h_\tau:\ssS\to\mathbb R_+\setminus\{0\}$ such that for each $t\in T$,
\equh\label{eq:fttau}
f_{t+\tau}(s) = h_\tau(s)\left(f_t\circ\Phi_\tau\right)(s)\,, \quad \ssL\ae
\eque
and
\equh\label{eq:dmu}
\dfrac\ssL{\Phi_\tau}\ssL(s) = h_\tau(s)^\alpha,\quad \ssL\ae\,.
\eque
Since, for every $t$, $\tau_1,\tau_2\in T$, we have two ways expressing $f_{t+\tau_1+\tau_2}$:
\[
f_{t+\tau_1+\tau_2}  =  f_{(t+\tau_1)+\tau_2} = (h_{\tau_2})(f_{t+\tau_1}\circ\Phi_{\tau_2}) 
 =  (h_{\tau_2})(h_{\tau_1}\circ\Phi_{\tau_2})(f_t\circ\Phi_{\tau_1}\circ\Phi_{\tau_2})\,, \ssL\ae
\]
and
\[
f_{t+\tau_1+\tau_2} = (h_{\tau_1+\tau_2})(f_t\circ\Phi_{\tau_1+\tau_2})\,,\ssL\ae\,,
\]
it follows, by the uniqueness of $\Phi_\tau$ and $h_\tau$, that for every $\tau_1,\tau_2\in T$,
\equh
h_{\tau_1+\tau_2} = (h_{\tau_2})(h_{\tau_1}\circ\Phi_{\tau_2})\,,\ssL\ae\,,
\eque
and
\equh\label{eq:phitau}
\Phi_{\tau_1+\tau_2} = \Phi_{\tau_1}\circ\Phi_{\tau_2}\,,\ssL\ae\,.
\eque
To complete the proof, we will establish a modification $\phi$ of $\Phi$ such that $\phi$ is measurable on $T\times\ssS$ and
\[
\Phi_t(s) = \phi(t,s)\,,\ssL\ae\,,\forall t\in T\,.
\]
If $T = \mathbb Z$, then one can modify $\indt\Phi$ to have~\eqref{eq:phitau} hold everywhere for all $\tau_1,\tau_2$, making $\indt\Phi$ a flow. When $T = \mathbb R$, 
by Theorem~1 in~\cite{mackey62point}, in order for $\indt\Phi$ to have a measurable version $\indt\phi$, it is enough to check that the map
\[
t\mapsto\tilde\nu\left(\left[\Phi_t\inv(B)\right]\right)
\]
is measurable for every finite measure $\tilde\nu$ on $\calB_\ssL$ (the measure algebra induced by $(\calB_\ssS,\ssL)$). It is clear that $\tilde\nu$ defines a finite measure $\nu$ on $\ssB$ such that $\nu(B) = \tilde\nu([B])$ and we have $\nu\ll\ssL$. Put $k = d\nu/d\ssL$. It is equivalent to show that
\equh\label{eq:tmapstoint}
t\mapsto\int_\ssS\ind_B\left(\Phi_t(s)\right)k(s)\ssL(ds)
\eque
is measurable for each $B\in\ssB$. Indeed, it is enough to show that $(t,s)\mapsto\ind_B(\Phi_t(s))$ is a measurable function of $(t,s)$ for each $B\in\ssB$. Choose a function $g = \ssF$ defined in~\eqref{eq:ssf} and $g_n = \bigvee_{i=1}^nc_{ni}f_{t_{ni}}\in\msspan\{f_t,t\in T\}$, such that $g_n\to g\,,\ssL\ae$. In view of~\eqref{rep:flowRepSS1}, for each $\tau\in T$, 
\[
h_\tau(s)g_n\circ\Phi_\tau(s) = \bigvee_{i=1}^nc_{ni}f_{t_{ni}+\tau}(s)\,,\ssL\ae\ \  s\in \ssS\,.
\]
Observe that the r.h.s. is a measurable function of $(\tau,s)$ for each $n\in\mathbb N$ and the l.h.s. converges $\ssL\ae$ as $n\to\infty$, for all $t\in T$. It follows that there exists a measurable function $(\tau,s)\mapsto g_\tau(s)$ such that, for each $\tau\in T$,
\equh\label{eq:htaus}
h_\tau(s)g\circ\Phi_\tau(s) = g_\tau(s)\,,\ssL\ae\,.
\eque
Now, observe that since $\indt f$ is minimal, for every $B\in\ssB$ there exist $\inftydots t \in T$ and $A\in\mathbb R^{\mathbb N}$ such that $B = \{s:(f_{t_1}(s)/g(s),f_{t_2}(s)/g(s),\dots)\in A\}\mod\ssL$. Note that~\eqref{eq:fttau} and~\eqref{eq:htaus} imply
\[
\frac{f_t\circ\Phi_\tau(s)}{g\circ\Phi_\tau(s)} = \frac{f_{t+\tau}(s)}{g_{\tau}(s)}\,,\ \ssL\ae\,.
\]
It follows that
\[
\ind_B(\Phi_\tau(s)) = \ind_A\left(f_{t_1+\tau}(s)/g_\tau(s),f_{t_2+\tau}(s)/g_\tau(s),\dots\right)\,,\ssL\ae\ \ s\in \ssS\,.
\]
We have thus shown that the map in~\eqref{eq:tmapstoint} is measurable.\ifthenelse{\boolean{qedTrue}}{\qed}{}
\end{proof}
\begin{proof}[Proof of Proposition~\ref{prop:equiFlow}] 
\itemnumber i The fact that $\indt{f\topp2}$ is another spectral representation of $\indt X$ can be verified by checking
\[
\left\|\bigvee c_jf_{t_j}\topp2\right\|_{L^\alpha_+(S_2,\mu_2)} = \left\| \bigvee c_jf_{t_j}\topp1 \right\|_{L^\alpha_+(S_1,\mu_1)}\,.
\]
\itemnumber {ii} By Corollary~\ref{coro:uniquePointMapping}, there exists measurable and invertible point mapping $\Phi:S_2\to S_1$ such that we have two different ways relating $f_{t+\tau}\topp2$ and $f_t\topp1$:
\[
f_{t+\tau}\topp2 = \Big(\dfrac{\mu_2}{\phi_\tau\topp2}{\mu_2}\Big)^\alpha f_t\topp2\circ\phi_\tau\topp2 = \Big(\dfrac{\mu_1}{\Phi\circ\phi_\tau\topp2}{\mu_2}\Big)^\alpha f_t\topp1\circ\Phi\circ\phi_\tau\topp2\,,\mu_2\ae\,,
\]
and
\[
f_{t+\tau}\topp2 = \Big(\dfrac{\mu_1}{\Phi}{\mu_2}\Big)^\alpha f_{t+\tau}\topp2\circ\phi_\tau\topp2 = \Big(\dfrac{\mu_1}{\phi_\tau\topp1\circ\Phi}{\mu_2}\Big)^\alpha f_t\topp1\circ\phi_\tau\topp1\circ\Phi\,,\mu_2\ae.
\]
By the uniqueness of determining flow, we have 
$
\phi_\tau\topp1\circ\Phi = \Phi\circ\phi_\tau\topp2\,,\mu_2\ae.
$
\end{proof}

\bibliographystyle{abbrv}

\begin{thebibliography}{10}

\bibitem{aaronson97introduction}
J.~Aaronson.
\newblock {\em An Introduction to Infinite Ergodic Theory}.
\newblock American Mathematical Society, 1997.

\bibitem{balkema77max}
A.~A. Balkema and S.~I. Resnick.
\newblock Max-infinite divisibility.
\newblock {\em Journal of Applied Probability}, 14(2):309--319, 1977.

\bibitem{brown77extreme}
B.~M. Brown and S.~I. Resnick.
\newblock Extreme values of independent stochastic processes.
\newblock {\em Journal of Applied Probability}, 14(4):732--739, 1977.

\bibitem{cohn72measurable}
D.~L. Cohn.
\newblock Measurable choice of limit points and the existence of separable and
  measurable processes.
\newblock {\em Zeitschrift f{\"u}r Wahrscheinlichkeitstheorie und verwandte
  Gebiete}, 22:161--165, 1972.

\bibitem{davis93prediction}
R.~A. Davis and S.~I. Resnick.
\newblock Prediction of stationary max-stable processes.
\newblock {\em Ann. Appl. Probab.}, 3(2):497--525, 1993.

\bibitem{dehaan78characterization}
L.~de~Haan.
\newblock A characterization of multidimensional extreme--value distributuions.
\newblock {\em Sankhy\-a (Statistics). The Indian Journal of Statistics. Series
  A}, 40(1):85--88, 1978.

\bibitem{dehaan84spectral}
L.~de~Haan.
\newblock A spectral representation for max-stable processes.
\newblock {\em Ann. Probab.}, 12(4):1194--1204, 1984.

\bibitem{dehaan86stationary}
L.~de~Haan and J.~Pickands~III.
\newblock {Stationary min-stable stochastic processes.}
\newblock {\em Probab. Theory Relat. Fields}, 72:477--492, 1986.

\bibitem{doob53stochastic}
J.~L. Doob.
\newblock {\em Stochastic Processes}.
\newblock Wiley, New York, 1953.

\bibitem{gine90max}
E.~Gin\'e, M.~G.~Hahn, and P.~Vatan.
\newblock Max-infinitely divisible and max-stable sample continuous processes.
\newblock {\em Probability Theory and Related Fields}, 87(2):139--165, 1990.

\bibitem{halmos50measure}
P.~R. Halmos.
\newblock {\em Measure Theory}.
\newblock Van Nostrand, Princeton, NJ, 1950.

\bibitem{hardin81isometries}
C.~D. Hardin, Jr.
\newblock Isometries on subspaces of $l^p$.
\newblock {\em Indiana University Mathematics Journal}, 30:449--465, 1981.

\bibitem{hardin82spectral}
C.~D. Hardin, Jr.
\newblock On the spectral representation of symmetric stable processes.
\newblock {\em Journal of Multivariate Analysis}, 12:385--401, 1982.

\bibitem{hida93gaussian}
T.~Hida and M.~Hitsuda.
\newblock {\em Gaussian processes}, volume 120 of {\em Translations of
  Mathematical Monographs}.
\newblock American Mathematical Society, Providence, RI, 1993.
\newblock Translated from the 1976 Japanese original by the authors.

\bibitem{kabluchko08spectral}
Z.~Kabluchko.
\newblock Spectral representations of sum-- and max--stable processes.
\newblock {\em submitted}, 2008.

\bibitem{kabluchko08stationary}
Z.~Kabluchko, M.~Schlather, and L.~de~Haan.
\newblock Stationary max--stable fields associated to negative definite
  functions.
\newblock {\em submitted}, 2008.

\bibitem{krengel69}
U.~Krengel.
\newblock Darstellungss\"atze f\"ur {S}tr\"omungen und {H}albstr\"omungen.
  {II}.
\newblock {\em Math. Ann.}, 182:1--39, 1969.

\bibitem{Krengel85ergodic}
U.~Krengel.
\newblock {\em Ergodic Theorems}.
\newblock de Gruyter, Berlin, 1985.

\bibitem{lamperti58isometries}
J.~Lamperti.
\newblock On the isometries of certain function-spaces.
\newblock {\em Pacific Journal of Mathematics}, 8, 1958.

\bibitem{mackey62point}
G.~W. Mackey.
\newblock Point realizations of transformation groups.
\newblock {\em Illinois Journal of Mathematics}, 6:327--335, 1962.

\bibitem{oodaira72strassen}
H.~Oodaira.
\newblock On {S}trassen's version of the law of the iterated logarithm for
  {G}aussian processes.
\newblock {\em Z. Wahrscheinlichkeitstheorie und Verw. Gebiete}, 21:289--299,
  1972.

\bibitem{pipiras02structure}
V.~Pipiras and M.~S. Taqqu.
\newblock The structure of self--similar stable mixed moving averages.
\newblock {\em Ann. Probab.}, 30, 2002.

\bibitem{pipiras04stable}
V.~Pipiras and M.~S. Taqqu.
\newblock Stable stationary processes related to cyclic flows.
\newblock {\em Ann. Probab.}, 32(3A):2222--2260, 2004.

\bibitem{resnick87extreme}
S.~I. Resnick.
\newblock {\em Extreme Values, Regular Variation and Point Processes}.
\newblock Springer-Verlag, New York, 1987.

\bibitem{resnick91random}
S.~I. Resnick and R.~Roy.
\newblock Random usc functions, max-stable processes and continuous choice.
\newblock {\em Ann. Appl. Probab.}, 1(2):267--292, 1991.

\bibitem{rosinski94uniqueness}
J.~Rosi\'nski.
\newblock On the uniqueness of spectral representation of stable processes.
\newblock {\em Journal of Theoretical Probability}, 7:615--634, 1994.

\bibitem{rosinski95structure}
J.~Rosi\'nski.
\newblock On the structure of stationary stable processes.
\newblock {\em Ann. Probab.}, 23(3):1163--1187, 1995.

\bibitem{rosinski00decomposition}
J.~Rosi\'nski.
\newblock Decompostion of stationary $\alpha$--stable random fields.
\newblock {\em Ann. Probab.}, 28:1797--1813, 2000.

\bibitem{rosinski06minimal}
J.~Rosi\'nski.
\newblock Minimal integral representations of stable processes.
\newblock {\em Probability and Mathematical Statistics}, 26:121--142, 2006.

\bibitem{rosinski96classes}
J.~Rosi{\'n}ski and G.~Samorodnitsky.
\newblock Classes of mixing stable processes.
\newblock {\em Bernoulli}, 2(4):365--377, 1996.

\bibitem{samorodnitsky05null}
G.~Samorodnitsky.
\newblock Null flows, positive flows and the structure of stationary symmetric
  stable processes.
\newblock {\em Ann. Probab.}, 33:1782--1803, 2005.

\bibitem{samorodnitsky94stable}
G.~Samorodnitsky and M.~S. Taqqu.
\newblock {\em Stable Non-Gaussian Random Processes}.
\newblock Chapman \& Hall, 1994.

\bibitem{sikorski64boolean}
R.~Sikorski.
\newblock {\em Boolean Algebras}.
\newblock Academic Press, New York, 1964.

\bibitem{stoev08ergodicity}
S.~A. Stoev.
\newblock On the ergodicity and mixing of max-stable processes.
\newblock {\em Stochastic Process. Appl.}, 118(9):1679--1705, 2008.

\bibitem{stoev06extremal}
S.~A. Stoev and M.~S. Taqqu.
\newblock Extremal stochastic integrals: a parallel between max-stable and
  alpha-stable processes.
\newblock {\em Extremes}, 8(3):237--266, 2006.

\end{thebibliography}

\end{document}